\documentclass[leqno]{siamart190516}

\usepackage{lipsum}
\usepackage{amsfonts}
\usepackage{graphicx}
\usepackage{epstopdf}
\ifpdf
  \DeclareGraphicsExtensions{.eps,.pdf,.png,.jpg}
\else
  \DeclareGraphicsExtensions{.eps}
\fi

\usepackage{tikz}
\usetikzlibrary{3d}
\usetikzlibrary{patterns}
\usetikzlibrary{calc}
\usetikzlibrary{arrows}

\usetikzlibrary{matrix}
\usetikzlibrary{positioning}

\usepackage{pgfplots}
\usepackage{pgfplotstable}
\usepackage{etoolbox}
\usetikzlibrary{decorations.pathreplacing}

\definecolor{mpiblue}{HTML}{33a5c3}

\usepackage{amsmath}
\pgfplotsset{compat=newest}

\usepackage{xcolor}
\usepackage{amssymb}
\usepackage{stmaryrd}
\usepackage{booktabs}
\usepackage{adjustbox}
\usepackage{caption}
\usepackage{subcaption}
\usepackage{mathtools}

\newcommand{\rulesep}{\unskip\ \vrule\ }

\newsiamremark{remark}{Remark}
\newsiamremark{hypothesis}{Hypothesis}
\crefname{hypothesis}{Hypothesis}{Hypotheses}
\newsiamthm{claim}{Claim}

\headers{Parallel Algorithms for TT Arithmetic}{H. Al Daas, G. Ballard, P. Benner}

\title{Parallel Algorithms for Tensor Train Arithmetic\thanks{Submitted to the editors - -, 2020.}}

\author{Hussam Al Daas\thanks{Department of Computational Methods in Systems and Control Theory, Max Planck Institute for Dynamics of Complex Technical Systems, Magdeburg, Germany
  (\email{aldaas@mpi-magdeburg.mpg.de}, \email{benner@mpi-magdeburg.mpg.de}).}
\and Grey Ballard\thanks{Computer Science Department, Wake Forest University, Winston Salem, North Carolina, USA
	(\email{ballard@wfu.edu}).}
\and Peter Benner\footnotemark[2]}

\usepackage{amsopn}

\newcommand{\datafile}{}
\newcommand{\captionsub}[1]{\caption{\scriptsize{#1}}}

\newcommand{\change}[1]{#1} 

\usepackage[mathscr]{eucal}

\usepackage{stmaryrd}

\usepackage{amsbsy}

\usepackage{bm}

\usepackage{braket}

\newcommand{\Tra}{\top}

\newcommand{\Real}{\mathbb{R}}

\newcommand{\V}[2][]{{\bm{#1\mathbf{\MakeLowercase{#2}}}}} 
 
\newcommand{\M}[2][]{{\bm{#1\mathbf{\MakeUppercase{#2}}}}} 
 
\newcommand{\T}[2][]{\boldsymbol{#1\mathscr{\MakeUppercase{#2}}}} 

\newcommand{\Mz}[3][]{\M[#1]{#2}_{(#3)}}

\newcommand{\Hada}{\ast} 
 
\newcommand{\Tentry}[1]{\ensuremath\T{#1}(i_1,\dots,i_N)}

\newcommand{\TT}[2]{\T{T}_{\T{#1},#2}} 
\newcommand{\slice}[2]{} 
\newcommand{\TTslice}[2]{ 
\ifx1#2  
	\renewcommand{\slice}[2]{\TT{#1}{1}(i_1,:)}
\else \ifx#2N 
	\renewcommand{\slice}[2]{\TT{#1}{N}(:,i_N)}
\else 
	\renewcommand{\slice}[2]{\TT{#1}{#2}(:,i_{#2},:)}
\fi \fi
\slice{#1}{#2} 
} 

\newcommand{\HOp}{\mathcal{H}}
\newcommand{\VOp}{\mathcal{V}}

\newcommand{\R}{\mathbb{R}} 

\pgfmathsetmacro{\none}{4.5}
\pgfmathsetmacro{\ntwo}{4.25}
\pgfmathsetmacro{\nthree}{5}
\pgfmathsetmacro{\nfour}{3}
\pgfmathsetmacro{\nfive}{4}
\pgfmathsetmacro{\rtwo}{1}
\pgfmathsetmacro{\rthree}{1.1}
\pgfmathsetmacro{\rfour}{1.05}
\pgfmathsetmacro{\rfive}{1.15}
\pgfmathsetmacro{\gap}{1}
\pgfmathsetmacro{\offset}{.25}
\pgfmathsetmacro{\vecwidth}{.1}

\newcommand{\nthcoredimone}{R_{n{-}2}}
\newcommand{\nthcoredimthree}{R_{n}}
\newcommand{\nextcoredimone}{R_{n}}
\newcommand{\nextcoredimthree}{R_{n{+}1}}
\newcommand{\nthcoredimtwo}{I_{n}}
\newcommand{\nthcoredimthreered}{L_{n}}

\newcommand{\drawTT}{
\draw (0,0,0) -- ++(-\rtwo,0,0) -- ++(0,-\none,0) -- ++(\rtwo,0,0) -- cycle;
\node at (-\rtwo-\offset,-.5*\none,0) {$I_1$};
\node at (-.5*\rtwo,\offset,0) {\small $R_1$};
\begin{scope}[shift={(\rthree+\gap,0,0)}]
	\draw (0,0,0) -- ++(-\rthree,0,0) -- ++(0,-\ntwo,0) -- ++(\rthree,0,0) -- cycle;
	\draw (0,0,0) -- ++(0,0,-\rtwo) -- ++(0,-\ntwo,0) -- ++(0,0,\rtwo) -- cycle;
	\draw (0,0,0) -- ++(-\rthree,0,0) -- ++(0,0,-\rtwo) -- ++(\rthree,0,0) -- cycle;
	\node at (-\rthree-\offset,-.5*\ntwo,0) {$I_2$};
	\node at (-\rthree-\offset,\offset,-.5*\rtwo) {\small $R_1$};
	\node at (-.5*\rthree,\offset,-\rtwo) {\small $R_2$};
\end{scope}
\begin{scope}[shift={(\rthree+\rfour+2*\gap,0,0)}]
	\draw (0,0,0) -- ++(-\rfour,0,0) -- ++(0,-\nthree,0) -- ++(\rfour,0,0) -- cycle;
	\draw (0,0,0) -- ++(0,0,-\rthree) -- ++(0,-\nthree,0) -- ++(0,0,\rthree) -- cycle;
	\draw (0,0,0) -- ++(-\rfour,0,0) -- ++(0,0,-\rthree) -- ++(\rfour,0,0) -- cycle;
	\node at (-\rfour-\offset,-.5*\nthree,0) {$I_3$};
	\node at (-\rfour-\offset,\offset,-.5*\rthree) {\small $R_2$};
	\node at (-.5*\rfour,\offset,-\rthree) {\small $R_3$};
\end{scope}
\begin{scope}[shift={(\rthree+\rfour+\rfive+3*\gap,0,0)}]
	\draw (0,0,0) -- ++(-\rfive,0,0) -- ++(0,-\nfour,0) -- ++(\rfive,0,0) -- cycle;
	\draw (0,0,0) -- ++(0,0,-\rfour) -- ++(0,-\nfour,0) -- ++(0,0,\rfour) -- cycle;
	\draw (0,0,0) -- ++(-\rfive,0,0) -- ++(0,0,-\rfour) -- ++(\rfive,0,0) -- cycle;
	\node at (-\rfive-\offset,-.5*\nfour,0) {$I_4$};
	\node at (-\rfive-\offset,\offset,-.5*\rfour) {\small $R_3$};
	\node at (-.5*\rfive,\offset,-\rfour) {\small $R_4$};
\end{scope}
\begin{scope}[shift={(\rthree+\rfour+\rfive+4*\gap,0,0)}]
	\draw (0,0,0) -- ++(0,0,-\rfive) -- ++(0,-\nfive,0) -- ++(0,0,\rfive) -- cycle;
	\node at (-\offset,-.5*\nfive,0) {$I_5$};
	\node at (-\offset,\offset,-.5*\rfive) {\small $R_4$};
\end{scope}
}

\newcommand{\drawHUnfolding}{
  \draw (0,0) rectangle (\rtwo, \rfive);
  \node at (0.5*\rtwo,-1.5*\offset) {\small $\nthcoredimthree$};
  \node at (-2.2*\offset,.5*\rfive) {\small $\nthcoredimone$};
  \begin{scope}[shift={(\rtwo,0)}]
    \draw (0,0) rectangle (1.5*\rtwo, \rfive);
    \node at (0.75*\rtwo,.5*\rfive) {$\cdots$};
  \end{scope}
  \begin{scope}[shift={(2.5\rtwo,0)}]
    \shade[top color=mpiblue!50,bottom color=mpiblue!50] (0,0) rectangle (\rtwo,\rfive);
    \draw (0,0) rectangle (\rtwo,\rfive);
    \node at (0.5*\rtwo,-1.5*\offset) {\small $\nthcoredimthree$};
  \end{scope}
  \begin{scope}[shift={(3.5*\rtwo,0)}]
    \draw (0,0) rectangle (1.5*\rtwo, \rfive);
    \node at (.75*\rtwo,.5*\rfive) {$\cdots$};
  \end{scope}
  \begin{scope}[shift={(5*\rtwo,0)}]
    \draw (0,0) rectangle (\rtwo, \rfive);
    \node at (0.5*\rtwo,-1.5*\offset) {\small $\nthcoredimthree$};
  \end{scope}
  \draw [thick, decoration = { brace, mirror, raise = 0.5cm}, decorate] (0,0) -- (6*\rtwo,0);
  \node at (3*\rtwo, -4.8*\offset) {\small $\nthcoredimtwo$};
}

\newcommand{\drawHUnfoldingO}{
  \draw (0,0) rectangle (\rthree, \rtwo);
  \node at (0.5*\rthree,-1.5*\offset) {\small $\nextcoredimthree$};
  \node at (6*\rthree+2.4*\offset,.5*\rtwo) {\small $\nextcoredimone$};
  \begin{scope}[shift={(\rthree,0)}]
    \draw (0,0) rectangle (1.5*\rthree, \rtwo);
    \node at (0.75*\rthree,.5*\rtwo) {$\cdots$};
  \end{scope}
  \begin{scope}[shift={(2.5*\rthree,0)}]
    \draw (0,0) rectangle (\rthree,\rtwo);
    \node at (0.5*\rthree,-1.5*\offset) {\small $\nextcoredimthree$};
  \end{scope}
  \begin{scope}[shift={(3.5*\rthree,0)}]
    \draw (0,0) rectangle (1.5*\rthree, \rtwo);
    \node at (.75*\rthree,.5*\rtwo) {$\cdots$};
  \end{scope}
  \begin{scope}[shift={(5*\rthree,0)}]
    \draw (0,0) rectangle (\rthree, \rtwo);
    \node at (0.5*\rthree,-1.5*\offset) {\small $\nextcoredimthree$};
  \end{scope}
}

\newcommand{\drawHUnfoldingOR}{
  \draw (0,0) rectangle (\rthree, \rtwo);
  \node at (0.5*\rthree,-1.5*\offset) {\small $\nextcoredimthree$};
  \node at (6*\rthree+2.4*\offset,.5*\rtwo) {\small $\nthcoredimthreered$};
  \begin{scope}[shift={(\rthree,0)}]
    \draw (0,0) rectangle (1.5*\rthree, \rtwo);
    \node at (0.75*\rthree,.5*\rtwo) {$\cdots$};
  \end{scope}
  \begin{scope}[shift={(2.5*\rthree,0)}]
    \draw (0,0) rectangle (\rthree,\rtwo);
    \node at (0.5*\rthree,-1.5*\offset) {\small $\nextcoredimthree$};
  \end{scope}
  \begin{scope}[shift={(3.5*\rthree,0)}]
    \draw (0,0) rectangle (1.5*\rthree, \rtwo);
    \node at (.75*\rthree,.5*\rtwo) {$\cdots$};
  \end{scope}
  \begin{scope}[shift={(5*\rthree,0)}]
    \draw (0,0) rectangle (\rthree, \rtwo);
    \node at (0.5*\rthree,-1.5*\offset) {\small $\nextcoredimthree$};
  \end{scope}
}

\newcommand{\drawVUnfolding}{
  \draw (0,0) rectangle (\rtwo, \rfive);
  \node at (0.5*\rtwo,-1.5*\offset) {\small $\nthcoredimthree$};
  \node at (\rtwo+2.5*\offset,.5*\rfive) {{\small $\nthcoredimone$}};
  \begin{scope}[shift={(0,\rfive)}]
    \draw (0,0) rectangle (\rtwo, 1.5*\rfive);
    \node at (0.5*\rtwo,.75*\rfive) {$\vdots$};
  \end{scope}
  \begin{scope}[shift={(0,2.5*\rfive)}]
    \shade[top color=mpiblue!50,bottom color=mpiblue!50] (0,0) rectangle (\rtwo,\rfive);
    \draw (0,0) rectangle (\rtwo,\rfive);
    \node at (\rtwo + 2.5*\offset,.5*\rfive) {{\small $\nthcoredimone$}};
  \end{scope}
  \begin{scope}[shift={(0,3.5*\rfive)}]
    \draw (0,0) rectangle (\rtwo, 1.5*\rfive);
    \node at (0.5*\rtwo,.75*\rfive) {$\vdots$};
  \end{scope}
  \begin{scope}[shift={(0,5*\rfive)}]
    \draw (0,0) rectangle (\rtwo, \rfive);
    \node at (\rtwo+2.5*\offset,.5*\rfive) {{\small $\nthcoredimone$}};
  \end{scope}
  \draw [thick, decoration = { brace, mirror, raise = 0.5cm}, decorate] (\rtwo + 2*\offset,0) -- (\rtwo + 2*\offset,6*\rfive);
  \node at (\rtwo+7*\offset, 3*\rfive) {\small $\nthcoredimtwo$};
}

\newcommand{\drawVUnfoldingO}{
  \draw (0,0) rectangle (\rtwo, \rfive);
  \node at (0.5*\rtwo,-1.5*\offset) {\small $\nthcoredimthree$};
  \node at (-\rtwo+\offset,.5*\rfive) {{\small $\nthcoredimone$}};
  \begin{scope}[shift={(0,\rfive)}]
    \draw (0,0) rectangle (\rtwo, 1.5*\rfive);
    \node at (0.5*\rtwo,.75*\rfive) {$\vdots$};
  \end{scope}
  \begin{scope}[shift={(0,2.5*\rfive)}]
    \draw (0,0) rectangle (\rtwo,\rfive);
    \node at (-\rtwo+\offset,.5*\rfive) {{\small $\nthcoredimone$}};
  \end{scope}
  \begin{scope}[shift={(0,3.5*\rfive)}]
    \draw (0,0) rectangle (\rtwo, 1.5*\rfive);
    \node at (0.5*\rtwo,.75*\rfive) {$\vdots$};
  \end{scope}
  \begin{scope}[shift={(0,5*\rfive)}]
    \draw (0,0) rectangle (\rtwo, \rfive);
    \node at (-\rtwo+\offset,.5*\rfive) {{\small $\nthcoredimone$}};
  \end{scope}
}

\newcommand{\drawVUnfoldingOR}{
  \draw (0,0) rectangle (\rtwo, \rfive);
  \node at (0.5*\rtwo,-1.5*\offset) {\small $\nthcoredimthreered$};
  \node at (-\rtwo-\offset,.5*\rfive) {{\small $\nthcoredimone$}};
  \begin{scope}[shift={(0,\rfive)}]
    \draw (0,0) rectangle (\rtwo, 1.5*\rfive);
    \node at (0.5*\rtwo,.75*\rfive) {$\vdots$};
  \end{scope}
  \begin{scope}[shift={(0,2.5*\rfive)}]
    \draw (0,0) rectangle (\rtwo,\rfive);
    \node at (-\rtwo-\offset,.5*\rfive) {{\small $\nthcoredimone$}};
  \end{scope}
  \begin{scope}[shift={(0,3.5*\rfive)}]
    \draw (0,0) rectangle (\rtwo, 1.5*\rfive);
    \node at (0.5*\rtwo,.75*\rfive) {$\vdots$};
  \end{scope}
  \begin{scope}[shift={(0,5*\rfive)}]
    \draw (0,0) rectangle (\rtwo, \rfive);
    \node at (-\rtwo-\offset,.5*\rfive) {{\small $\nthcoredimone$}};
  \end{scope}
}

\newcommand{\drawTower}{
	\draw (0,0,0) -- ++(-\rfive,0,0) -- ++(0,-\ntwo,0) -- ++(\rfive,0,0) -- cycle;
	\draw (0,0,0) -- ++(0,0,-\rtwo) -- ++(0,-\ntwo,0) -- ++(0,0,\rtwo) -- cycle;
	\draw (0,0,0) -- ++(-\rfive,0,0) -- ++(0,0,-\rtwo) -- ++(\rfive,0,0) -- cycle;
	\node at (-\rfive-\offset,-.5*\ntwo,0) {$\nthcoredimtwo$};
  \node at (-\rfive-2*\offset,\offset,-.5*\rtwo) {\small $\nthcoredimone$};
  \node at (-.5*\rfive,\offset,-\rtwo) {\small $\nthcoredimthree$};
}

\begin{document}
\maketitle
\begin{abstract}
We present efficient and scalable parallel algorithms for performing mathematical operations for low-rank tensors represented in the tensor train (TT) format.
We consider algorithms for addition, elementwise multiplication, computing norms and inner products, \change{orthonormal}ization, and rounding (rank truncation).
These are the kernel operations for applications such as iterative Krylov solvers that exploit the TT structure.
  The parallel algorithms are designed for distributed-memory computation, and we \change{propose} a data distribution and strategy that parallelizes computations for individual cores within the TT format.
We analyze the computation and communication costs of the proposed algorithms to show their scalability, and we present numerical experiments that demonstrate their efficiency on both shared-memory and distributed-memory parallel systems.
For example, we observe better single-core performance than the existing MATLAB TT-Toolbox in rounding a 2GB TT tensor, and our implementation achieves a $34\times$ speedup using all 40 cores of a single node.
We also show nearly linear parallel scaling on larger TT tensors up to over 10,000 cores for all mathematical operations.
\end{abstract}



\section{Introduction}

Multi-dimensional data, or tensors, appear in a variety of applications where numerical values represent multi-way relationships.
The Tensor Train (TT) format is a low-rank representation of a tensor that has been applied to solving problems in areas such as parameter-dependent PDEs, stochastic PDEs, molecular simulations, uncertainty quantification, data completion, and classification \cite{BenDOS16,BenDOS20,DolS17,HacK09, KreKNT15, KreT10, Ose11, SavDWK14}. 
As the number of dimensions or modes of a tensor becomes large, the total number of data elements grows exponentially fast, which is known as the curse of dimensionality \cite{HacK09}.
Fortunately, it can be shown in many cases that the tensors exhibit low-rank structure and can be represented or approximated by significantly fewer parameters.
Low-rank tensor approximations allow for storing the data implicitly and performing arithmetic operations in feasible time and space complexity, avoiding the curse of dimensionality.

In contrast to the matrix case where the singular value decomposition (SVD) provides optimal low-rank representations, there are more diverse possibilities for low-rank representations of tensors \cite{KolB09}.
Various representations have been proposed, such as CP \cite{CC70,Harshman70}, Tucker \cite{Tucker66}, quantized tensor train \cite{Kho11}, and hierarchical Tucker \cite{HacK09}, in addition to TT \cite{Ose11}, and each has been demonstrated to be most effective in certain applications. 
The TT format\change{, which is also known as the matrix product state (MPS) in the computational physics and chemistry communities,} consists of a sequence of TT cores, one for each tensor dimension, and each core is a 3-way tensor except for the first and last cores, which are matrices.
The primary advantages of TT are that (1) the number of parameters in the representation is linear, rather than exponential, in the number of modes and (2) the representation can be computed to satisfy a specified approximation error threshold in a numerically stable way.

As these low-rank tensor techniques have been applied to larger and larger data sets, efficient sequential and parallel implementations of algorithms for computing and manipulating these formats have also been developed.
Toolboxes and libraries in productivity-oriented languages such as MATLAB and Python \cite{TensorToolbox,Tensorly,TTToolbox,Tensorlab} are available for moderately sized data, and parallel algorithms implemented in performance-oriented languages exist for computation of decompositions such as CP \cite{EH+19-TR,SRSK15,LCPSV17} and Tucker \cite{AusBT16,BalKK20,KU16,SK17} and operations such as tensor contraction \cite{SolMHSD14}, allowing for scalability to much larger data and numbers of processors.
\change{While efficient computation of TT approximations of explicit tensors has attracted recent attention \cite{BC+20,GK21,LYB20,RTB21,WYWRD21}}, no such \change{high-performance parallel implementations exist for approximating tensors already in TT format}.
\change{
In condensed matter computations, several advances have been made in parallelizing the density matrix renormalization group (DMRG) algorithm, which computes the ground-state eigenvector in MPS/TT format \cite{KanDTG19,LevSC20,StoW13}.
The modes' dimensions in these applications are very small and the TT ranks can be very large.
In contrast, applications from parameter-dependent PDEs, stochastic PDEs, uncertainty quantification, and molecular simulations \cite{BenDOS20,BenGW15,KreKNT15} yield computations with TT tensors having certain modes with very large dimensions and relatively small TT ranks.}
The goal of this work is to establish efficient and scalable algorithms for implementing the key mathematical operations on TT tensors 
\change{ for applications where at least one mode has a very large dimension and the TT ranks are relatively small to allow researchers to scale their models beyond the time and memory constraints when using current MATLAB and Python implementations.
}

We consider mathematical operations such as addition, Hadamard (elementwise) multiplication, computing norms and inner products, left- and right-\change{orthonormal}ization, as well as rounding (rank truncation).
These are the operations required to, for example, solve a structured linear system whose solution can be approximated well by a tensor in TT format using a Krylov method \cite{KreT10}.
As we will see in \Cref{sec:background}, mathematical operations can increase the \change{ranks of the TT representation of the result tensor}, which can then be recompressed, or rounded back to smaller ranks, in order to maintain feasible time and space complexity with some controllable loss of accuracy.
As a result, the rounding procedure (and the \change{orthonormal}ization it requires) is of prime importance in developing efficient and scalable TT algorithms.
We will assume throughout that full tensors are never formed explicitly.

In order to develop scalable parallel algorithms, we \change{propose} a data distribution and parallelization techniques that maintain computational load balance and attempt to minimize interprocessor communication, which is the most expensive operation on parallel machines in terms of both time and energy consumption.
As discussed in \Cref{sec:parallel-tt}, we distribute the slices of each TT core across all processors, where slices are matrices (or vectors) whose dimensions are determined by the low ranks of the TT representation.
This distribution allows for full parallelization of each core-wise computation and avoids the need for communication within slice-wise computations.
The \change{orthonormal}ization and rounding algorithms depend on parallel QR decompositions, and our approach enables the use of the Tall-Skinny QR algorithm, which is communication optimal for the matrix dimensions in this application \cite{DGHL12}.
We analyze the parallel computation and communication costs of each TT algorithm, demonstrating that the bulk of the computation is load balanced perfectly across processors.
The communication costs are independent of the original tensor dimensions, so their relative costs diminish with small ranks.

We verify the theoretical analysis and benchmark our C/MPI implementation on up to 256 nodes (10,240 cores) of a distributed-memory parallel platform in \Cref{sec:numerical_experiments}.
Our experiments are performed on synthetic data using tensor dimensions and ranks that arise in a variety of scientific and data analysis applications.
On a shared-memory system (one node of the system), we compare our TT-rounding implementation against the TT-Toolbox \cite{TTToolbox} in MATLAB and show that our implementation is 70\% more efficient using a single core and achieves up to a $34\times$ parallel speedup using all 40 cores on the node.
We also present strong scaling performance experiments for computing inner products, norms, \change{orthonormal}ization, and rounding using up to over 10K MPI processes.
The experimental results show that the time remains dominated by local computation even at that scale, allowing for nearly linear scaling for multiple operations, achieving for example a $97\times$ speedup of TT-rounding when scaling from 1 node to 128 nodes on a TT tensor with a 28 GB memory footprint.
We conclude in \Cref{sec:conclusion} and discuss limitations of our approaches and perspectives for future improvements.

\section{Notation and background}
\label{sec:background}

In this section, we review the tensor train (TT) format and present a brief overview of the notation and computational kernels associated with it.
Tensors are denoted by boldface Euler script letters (e.g. $\T{X}$), and matrices are denoted by boldface block letters (e.g. $\M A$).
The number $I_n$ for $1\leq n \leq N$ is referred to as the mode size or mode dimension, and we use $i_n$ to index that dimension.
The order of a tensor is its number of modes, e.g., the order of $\T{X}$ is $N$.
The $n$th TT core (described below) of a tensor $\T{X}$ is denoted by $\TT{X}{n}$.
We use MATLAB-style notation to obtain elements or sub-tensors, where a solitary colon (:) refers to the entire range of a dimension.
For example $\T{X}(i,j,k)$ is a tensor entry, $\T{X}(i,:,:)$ is a tensor slice (a matrix in this case), and $\T{X}(:,j,k)$ is a tensor fiber (a vector).

The mode-$n$ ``modal'' unfolding (or matricization or flattening) of a tensor $\T{X} \in \R^{I_1 \times I_2\times I_3}$ is the matrix $\Mz{X}{n} \in \R^{I_n \times \frac{I}{I_n}}$, where $I = I_1 I_2 I_3$.
In this case, the columns of the modal unfolding are fibers in that mode.
The mode-$n$ product or tensor-times-matrix operation is denoted by $\times_n$ and is defined so that the mode-$n$ unfolding of $\T{X} \times_n \M{A}$ is $\M{A}\Mz{X}{n}$.
We refer to \cite{KolB09,RVL12} for more details.

\change{The norm of a tensor is defined so that
$\|\T{X}\|^2 = \sum_{i_1,\dots,i_N} \Tentry{X}^2,$
which generalizes the vector 2-norm and matrix Frobenius norm.}

\subsection{TT tensors}

\begin{figure}
	\label{fig:representation_tt}
\begin{center}
  \begin{adjustbox}{height=0.3\linewidth}

\begin{tikzpicture}[scale=.7,textnode/.style={scale=.5}]
\pgfmathsetmacro{\ione}{.25*\none}
\fill[mpiblue!50] (0,-\ione,0) rectangle (-\rtwo,-\ione-\vecwidth,0);
\begin{scope}[shift={(\rthree+\gap,0,0)}]
	\pgfmathsetmacro{\itwo}{.6*\ntwo}
	\fill[mpiblue!50] (0,-\itwo,0) rectangle (-\rthree,-\itwo-\vecwidth,0);
	\fill[mpiblue!50] (0,-\itwo,0) -- ++(0,0,-\rtwo) -- ++(0,-\vecwidth,0) -- ++(0,0,\rtwo) -- cycle;
	\fill[mpiblue!50] (0,-\itwo,0) -- ++(-\rthree,0,0) -- ++(0,0,-\rtwo) -- ++(\rthree,0,0) -- cycle;
\end{scope}
\begin{scope}[shift={(\rthree+\rfour+2*\gap,0,0)}]
	\pgfmathsetmacro{\ithree}{.15*\nthree}
	\fill[mpiblue!50] (0,-\ithree,0) rectangle (-\rfour,-\ithree-\vecwidth,0);
	\fill[mpiblue!50] (0,-\ithree,0) -- ++(0,0,-\rthree) -- ++(0,-\vecwidth,0) -- ++(0,0,\rthree) -- cycle;
	\fill[mpiblue!50] (0,-\ithree,0) -- ++(-\rfour,0,0) -- ++(0,0,-\rthree) -- ++(\rfour,0,0) -- cycle;
\end{scope}
\begin{scope}[shift={(\rthree+\rfour+\rfive+3*\gap,0,0)}]
	\pgfmathsetmacro{\ifour}{.85*\nfour}
	\fill[mpiblue!50] (0,-\ifour,0) rectangle (-\rfive,-\ifour-\vecwidth,0);
	\fill[mpiblue!50] (0,-\ifour,0) -- ++(0,0,-\rfour) -- ++(0,-\vecwidth,0) -- ++(0,0,\rfour) -- cycle;
	\fill[mpiblue!50] (0,-\ifour,0) -- ++(-\rfive,0,0) -- ++(0,0,-\rfour) -- ++(\rfive,0,0) -- cycle;
\end{scope}
\begin{scope}[shift={(\rthree+\rfour+\rfive+4*\gap,0,0)}]
	\pgfmathsetmacro{\ifive}{.5*\nfive}
	\fill[mpiblue!50] (0,-\ifive,0) -- ++(0,0,-\rfive) -- ++(0,-\vecwidth,0) -- ++(0,0,\rfive) -- cycle;
\end{scope}
\drawTT
\end{tikzpicture}
  \end{adjustbox}
\end{center}
        \caption{Order-5 TT tensor with a particular slice from each TT core highlighted.  The chain product of these slices produces a scalar element of the full tensor with indices corresponding to the slices.}
\end{figure}
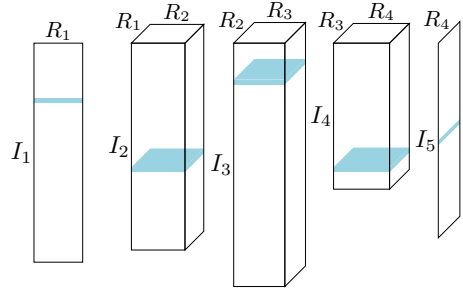

A tensor $\T{X} \in \R^{I_1 \times \cdots \times I_N}$ is in the TT format if there exist strictly positive integers $R_0, \ldots, R_N$ with $R_0 = R_N = 1$ and $N$ order-3 tensors $\TT{X}{1}, \ldots, \TT{X}{N} $, called TT cores, with $\TT{X}{n} \in \R^{R_{n-1} \times I_n \times R_n}$, such that:
\begin{displaymath}
	\Tentry{X} = \TTslice{X}{1} \cdots \TTslice{X}{n} \cdots \TTslice{X}{N}.
\end{displaymath}
We note that because $R_0=R_N=1$, the first and last TT cores are (order-2) matrices so $\TTslice{X}{1} \in \R^{R_1}$ and $\TTslice{X}{N} \in \R^{R_{N-1}}$.
The $R_{n-1} \times R_n$ matrix $\TTslice{X}{n}$ is referred to as the $i_n$th slice of the $n$th TT core of $\T{X}$, where $1\leq i_n\leq I_n$.
\Cref{fig:representation_tt} shows an illustration of an order-5 TT tensor.

Due to the multiplicative formulation of the TT format, the cores of a TT tensor are not unique.
For example, let $\T{X}$ be a TT tensor and $\M{M} \in \R^{R_{n} \times R_{n}}$ be an invertible matrix.
Then, the TT tensor $\T{Y}$ defined such that
\begin{equation*}
\Tentry{Y} = \TTslice{X}{1}\cdots (\TTslice{X}{n} \M{M})\cdot (\M{M}^{-1} \TTslice{X}{n+1}) \cdots \TTslice{X}{N}
\end{equation*}
is equal to $\T{X}$.
Another important remark is the following:
\begin{multline}
\label{eq:nonunique}
\TTslice{X}{1}\cdots (\TTslice{X}{n} \M{M})\cdot \TTslice{X}{n+1} \cdots \TTslice{X}{N} = \\
\TTslice{X}{1}\cdots \TTslice{X}{n} \cdot (\M{M} \TTslice{X}{n+1}) \cdots \TTslice{X}{N} 
\end{multline}
where $\M M$ in this case need not be invertible.
Thus, we can ``pass'' a matrix between adjacent cores without changing the tensor.
This property is used to \change{orthonormal}ize TT cores as we will see in \Cref{sec:orth}.

%

\subsection{Unfolding TT cores}

In order to express the arithmetic operations on TT cores using linear algebra, we will often use two specific matrix unfoldings of the 3D tensors.
The \emph{horizontal unfolding} of TT core $\TT{X}{n}$ corresponds to the concatenation of the slices $\TTslice{X}{n}$ for $i_n = 1,\ldots, I_n$ horizontally.
We denote the corresponding operator by $\HOp$, so that $\HOp(\TT{X}{n})$ is an $R_{n-1} \times R_n I_n$ matrix.
The \emph{vertical unfolding} corresponds to the concatenation of the slices $\TTslice{X}{n}$ for $i_n = 1,\ldots, I_n$ vertically.
We denote the corresponding operator by $\VOp$, so that $\VOp(\TT{X}{n})$ is an $R_{n-1} I_n \times R_n$ matrix.
These unfoldings are illustrated in \Cref{fig:unfolding}.

Note that the horizontal unfolding is equivalent to the modal unfolding with respect to the 1st mode, often denoted with subscript $(1)$ to denote the mode that corresponds to rows \cite{KolB09}.
Similarly, the vertical unfolding is the transpose of the modal unfolding with respect to the 3rd mode, which also corresponds to the more general unfolding that maps the first two modes to rows and the third mode to columns, denoted with subscript $(1{:}2)$ to denote the modes that correspond to rows \cite{PTC13a}.
These connections are important for the linearization of tensor entries in memory and our efficient use of BLAS and LAPACK, discussed in \Cref{sec:distribution}.

\renewcommand{\nthcoredimone}{R_{n-1}}
\renewcommand{\nthcoredimthree}{R_{n}}

\begin{figure}
  \centering
  \begin{adjustbox}{height=0.15\linewidth}

  \begin{tabular}{ccc}
    \begin{minipage}{0.23\linewidth}
      \begin{tikzpicture}[scale=.7,textnode/.style={scale=.5}]
        \pgfmathsetmacro{\ione}{.25*\none}
        \pgfmathsetmacro{\itwo}{.6*\ntwo}
        \fill[mpiblue!50] (0,-\itwo,0) rectangle (-\rfive,-\itwo-\vecwidth,0);
        \shade[top color=mpiblue!50,bottom color=mpiblue!50] (0,-\itwo,0) -- ++(0,0,-\rtwo) -- ++(0,-\vecwidth,0) -- ++(0,0,\rtwo) -- cycle;
        \shade[top color=mpiblue!50,bottom color=mpiblue!50] (0,-\itwo,0) -- ++(0,0,-\rtwo) -- ++(0,-\vecwidth,0) -- ++(0,0,\rtwo) -- cycle;
        \shade[top color=mpiblue!50,bottom color=mpiblue!50] (0,-\itwo,0) -- ++(-\rfive,0,0) -- ++(0,0,-\rtwo) -- ++(\rfive,0,0) -- cycle;
        \drawTower
      \end{tikzpicture}\\
      \centering{\footnotesize 
        $\TT{X}{n} \in \Real^{\nthcoredimone \times \nthcoredimtwo \times \nthcoredimthree}$ \\ 
        is a \emph{TT core} 
      }
    \end{minipage}&
    \begin{minipage}{0.42\linewidth}
      \begin{tikzpicture}[scale=.7,textnode/.style={scale=.5}]
        \drawHUnfolding
      \end{tikzpicture}
      \centering{\footnotesize 
        $\HOp(\TT{X}{n}) \in \Real^{\nthcoredimone \times \nthcoredimtwo\nthcoredimthree}$ \\  is \emph{horizontal unfolding}
      }
    \end{minipage}&
    \begin{minipage}{0.25\linewidth}
      \begin{tikzpicture}[scale=.7,textnode/.style={scale=.5}]
        \centering
        \drawVUnfolding
      \end{tikzpicture}
      \centering{\footnotesize 
        $\VOp(\TT{X}{n}) \in \Real^{\nthcoredimone\nthcoredimtwo \times \nthcoredimthree}$ \\ is \emph{vertical unfolding}
      }
    \end{minipage}
  \end{tabular}
  \end{adjustbox}
  \caption{Horizontal and vertical unfoldings of a TT core.}
  \label{fig:unfolding}
\end{figure}
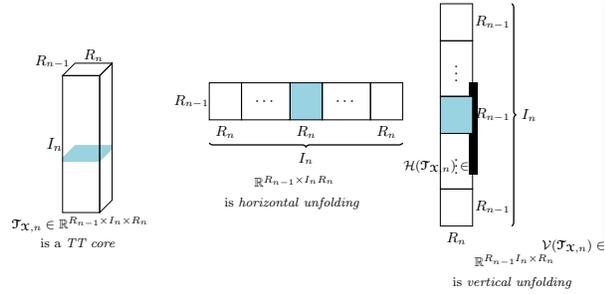

%

\subsection{TT \change{Orthonormal}ization}
\label{sec:orth}

Different types of \change{orthonormal}ization can be defined for TT tensors.
We focus in this paper on left and right \change{orthonormal}izations which are required in the rounding procedure.
We use the terms column and row \change{orthonormal} to refer to matrices that have \change{orthonormal} columns and \change{orthonormal} rows, respectively, so that a matrix $\M{Q}$ is column \change{orthonormal} if $\M{Q}^\Tra\M{Q}=\M{I}$ and row \change{orthonormal} if $\M{Q}\M{Q}^\Tra=\M{I}$.

A TT tensor is said to be right \change{orthonormal} if $\HOp(\TT{X}{n})$ is row \change{orthonormal} for $n=2,\ldots, N$ (all but the first core).
On the other hand, a tensor is said to be left \change{orthonormal} if $\VOp(\TT{X}{n})$ is column \change{orthonormal} for $n=1, \ldots, N-1$ (all but the last core).
More generally, we define a tensor to be $n$-right \change{orthonormal} if the horizontal unfoldings of cores $n+1, \ldots, N$ are all row \change{orthonormal}, and a tensor is $n$-left \change{orthonormal} if the vertical unfoldings of cores $1, \dots, n-1$ are all column \change{orthonormal}.

These definitions correspond to the fact that the tensor that represents the contraction of these sets of TT cores inherits their \change{orthonormal}ity.
For example, let $\T X$ be a right-\change{orthonormal} TT tensor, then we can write $\Mz{X}{1} = \TT{X}{1} \Mz{Z}{1}$, where $\T Z$ is a $R_1 \times I_2 \times \cdots \times I_N$ tensor whose entries are given by
$$\T{Z}(r_1,i_2,\dots,i_N) = \TT{X}{2}(r_1,i_2,:) \cdot \TTslice{X}{3} \cdots \TTslice{X}{n} \cdots \TTslice{X}{N}.$$
The 1st modal unfolding of $\T{Z}$ is row \change{orthonormal}, as shown below \cite[Lemma 3.1]{Ose11}:
\small
\change{\begin{align*}
\Mz{Z}{1}\Mz{Z}{1}^\Tra &= \sum_{i_2,\dots,i_N} \T{Z}(:,i_2,\dots,i_N) \T{Z}(:,i_2,\dots,i_N)^\Tra \\
 &= \sum_{i_2,\dots,i_N} \underbrace{\TTslice{X}{2} \cdots \TTslice{X}{N}}_{\change{\T{Z}(:,i_2,\dots,i_N)}}  \underbrace{\TTslice{X}{N}^\Tra \cdots \TTslice{X}{2}^\Tra}_{\T{Z}(:,i_2,\dots,i_N)^\Tra} \\
 &= \begin{aligned}[t] \sum_{i_2\change{,\dots,i_{N{-}1}}} \TTslice{X}{2} \cdots \TTslice{X}{N{-}1} \left(\sum_{i_N}  \TTslice{X}{N} \TTslice{X}{N}^\Tra \right) \\ \cdot \TTslice{X}{N{-}1}^\Tra \cdots \TTslice{X}{2}^\Tra \end{aligned} \\
 &= \begin{aligned}[t] \sum_{i_2,\dots,i_{N{-}1}} \TTslice{X}{2} \cdots \TTslice{X}{N{-}1} \underbrace{\HOp(\TT{X}{N}) \HOp(\TT{X}{N})^\Tra}_{I_{R_{N{-}1}}}  \\ \cdot \TTslice{X}{N{-}1}^\Tra \cdots \TTslice{X}{2}^\Tra \end{aligned} \\
 &= \sum_{i_2\change{,\dots,i_{N{-}1}}} \TTslice{X}{2} \cdots \TTslice{X}{N{-}1} \TTslice{X}{N{-}1}^\Tra \cdots \TTslice{X}{2}^\Tra \\
  &=\begin{aligned}[t] \sum_{i_2,\dots,i_{N{-}2}} \TTslice{X}{2} \cdots \TTslice{X}{N{-}2} \Bigg(\sum_{i_{N-1}}  \TTslice{X}{N{-}1} \\ \cdot \TTslice{X}{N{-}1}^\Tra \Bigg)
    \TTslice{X}{N{-}2}^\Tra \cdots \TTslice{X}{2}^\Tra \end{aligned} \\
 &= \begin{aligned}[t] \sum_{i_2,\dots,i_{N{-}2}} \TTslice{X}{2} \cdots \TTslice{X}{N{-}2} \underbrace{\HOp(\TT{X}{N-1}) \HOp(\TT{X}{N-1})^\Tra}_{I_{R_{N{-}2}}} \\\cdot \TTslice{X}{N{-}2}^\Tra \cdots \TTslice{X}{2}^\Tra \end{aligned} \\
  &= \sum_{i_2,\dots,i_{N{-}2}} \TTslice{X}{2} \cdots \TTslice{X}{N{-}2} \TTslice{X}{N{-}2}^\Tra  \cdots \TTslice{X}{2}^\Tra \\
 &= \cdots = I_{R_1}.
\end{align*}}
\normalsize
Similar arguments show that the 1st modal unfolding of the tensor representing the last $N-n$ cores of an $n$-right \change{orthonormal} TT tensor is row \change{orthonormal} and that the last modal unfolding of the tensor representing the first $n-1$ cores of an $n$-left \change{orthonormal} TT tensor is row \change{orthonormal}.

Given a TT tensor, we can \change{orthonormal}ize it by exploiting the non-uniqueness of TT tensors expressed in \Cref{eq:nonunique}.
That is, we can right- or left-\change{orthonormal}ize a TT core using a QR decomposition of one of its unfoldings and pass its triangular factor to its neighbor core without changing the represented tensor.
By starting from one end and repeating this process on each core in order, we can obtain a left or right \change{orthonormal} TT tensor, as shown in \Cref{alg:TT-rorthonormalization} (for right \change{orthonormal}ization).

\begin{algorithm}
\caption{TT-right-\change{orthonormal}ization}
\label{alg:TT-rorthonormalization}
\begin{algorithmic}[1]
	\Require{A TT tensor $\T{X}$}
	\Ensure{A right \change{orthonormal} TT tensor $\T{Y}$ equivalent to $\T{X}$}
	\Function{$\T{Y} =$ Right-\change{Orthonormal}ization}{$\T{X}$}
		\State Set $\TT{Y}{N} = \TT{X}{N}$
		\For{$n=N$ down to $2$}
			\State $[\HOp(\TT{Y}{n})^\Tra,\M{R}] = \textsc{QR}(\HOp(\TT{Y}{n})^\Tra)$ \Comment{QR factorization}
			\State $\VOp(\TT{Y}{n-1}) = \VOp(\TT{X}{n-1}) \M{R}^\Tra$ \Comment{$\TT{Y}{n-1} = \TT{X}{n-1} \times_3 \M{R}^\Tra$}
		\EndFor
	\EndFunction
\end{algorithmic}
\end{algorithm}


We note that the norm of a right- or left-\change{orthonormal} TT tensor can be cheaply computed, based on the idea that post-multiplication by a matrix with orthonormal rows or pre-multiplication by a matrix with orthonormal columns does not affect the Frobenius norm of a matrix.
Thus, we have that $\|\T X\| = \|\TT{X}{1}\|_F$ provided that $\Mz{Z}{1}$ has orthonormal rows, and $\|\T X\| = \|\TT{X}{N} \|_F$ if $\T X$ is left \change{orthonormal}.

\subsection{TT Rounding}
\label{sec:seq-tt-rounding}

\change{Orthonormal}ization plays an essential role in compressing the TT format of a tensor (decreasing the TT ranks $R_n$) \cite{Ose11}.
This compression is known as TT rounding and is given in \Cref{alg:TT-rounding}.

The intuition for rounding can be expressed in matrix notation as follows.
Suppose we have a matrix represented by a product 
\begin{equation}
\label{eq:quadprodmatrix}
\M{A} = \M{Q} \M{B} \M{C} \M{Z},
\end{equation}
where $\M{Q}$ and $\M{Z}$ are column and row \change{orthonormal}, respectively.
Then the truncated SVD of $\M{A}$ can be readily expressed in terms of the truncated SVD of $\M{BC}$. 
In our case, $\M{B}$ is tall and skinny and $\M{C}$ is short and wide, so the rank is bounded by their shared dimension.
To truncate the rank, one can row-\change{orthonormal}ize $\M{C}$ and then perform a truncated SVD of $\M{B}$ (or vice-versa).
That is, if we compute $\M{R}_C\M{Q}_C=\M{C}$ and $\M{U}_B\M{\Sigma}_B\M{V}_B^\Tra = \M{B}\M{R}_C$, then to round $\M{A}$ we can replace $\M{B}$ with $\M[\hat]{U}_B$ and $\M{C}$ with $\M[\hat]{\Sigma}_B \M[\hat]{V}_B^\Tra \M{Q}_C$, where $\M[\hat]{U}_B \M[\hat]{\Sigma}_B \M[\hat]{V}_B^\Tra$ is the SVD truncated to the desired tolerance.
 
\begin{figure}
\centering
\begin{adjustbox}{height=0.3\linewidth}

\pgfmathsetmacro{\mm}{7} 
\pgfmathsetmacro{\nn}{7.5} 
\pgfmathsetmacro{\im}{4} 
\pgfmathsetmacro{\in}{2.5} 
\pgfmathsetmacro{\rr}{.75} 
\pgfmathsetmacro{\ss}{.75} 
\pgfmathsetmacro{\vv}{.5} 

\begin{tikzpicture}[scale=.7,textnode/.style={scale=.5}]
	\begin{scope}[shift={(0,0)}]
		\draw (0,0) rectangle ++(\im/4,-\mm/4);
		\draw (\im/4,-\mm/4) rectangle ++(\im/4,-\mm/4);
		\draw[thick,loosely dotted] (1.05*\im/2,-1.05*\mm/2)-- ++(.8*\im/4,-.8*\mm/4);
		\draw (3*\im/4,-3*\mm/4) rectangle ++(\im/4,-\mm/4);
		\node[rotate=90] at (3*\im/4-\vv/2,-7*\mm/8) {\tiny $I_1{\cdots} I_{n{-}1}$};
		\node at (7*\im/8,-3*\mm/4+\vv/2) {\tiny $R_{n{-}1}$};
		\draw (0,0) rectangle ++(\im,-\mm);
		\node at (\im/2,-\mm-\vv) {$\M{I}_{I_n} \otimes \Mz{Q}{1:n-1}$};
		\node[rotate=90] at (-\vv,-\mm/2) {\small $I_1\cdots I_n$};
		\node at (\im/2,\vv) {\small $I_nR_{n-1}$};
	\end{scope}
	\begin{scope}[shift={(\im+\ss,0)}]
		\draw (0,0) rectangle ++(\rr,-\im);
		\node at (\rr/2,-\im-\vv) {$\VOp(\TT{X}{n})$};
		\node at (\rr/2,\vv) {\small $R_{n}$};
	\end{scope}
	\begin{scope}[shift={(\im+\rr+2*\ss,0)}]
		\draw (0,0) rectangle ++(\in,-\rr);
		\node at (\in/2,-\rr-\vv) {$\HOp(\TT{X}{n+1})$};
		\node at (\in/2,\vv) {\small $I_{n+1}R_{n+1}$};
	\end{scope}
	\begin{scope}[shift={(\im+\rr+\in+3*\ss,0)}]
		\draw (0,0) rectangle ++(\nn/4,-\in/4);
		\draw (\nn/4,-\in/4) rectangle ++(\nn/4,-\in/4);
		\draw[thick,loosely dotted] (1.05*\nn/2,-1.05*\in/2)-- ++(.8*\nn/4,-.8*\in/4);		
		\draw (3*\nn/4,-3*\in/4) rectangle ++(\nn/4,-\in/4);
		\node at (7*\nn/8,-3*\in/4+\vv/2) {\tiny $I_{n{+}2}{\cdots} I_N$};
		\node at (3*\nn/4-\vv,-7*\in/8) {\tiny $R_{n+1}$};
		\draw (0,0) rectangle ++(\nn,-\in);
		\node at (\nn/2,-\in-\vv) {$\M{I}_{I_{n+1}} \otimes \Mz{Z}{1}$};
		\node at (\nn/2,\vv) {\small $I_{n+1}\cdots I_N$};
	\end{scope}
\end{tikzpicture}
\end{adjustbox}
  \caption{Visualization of identity \cref{eq:quadprodTT} for $\Mz{X}{1:n}$.}
\label{fig:quadprod}
\end{figure}
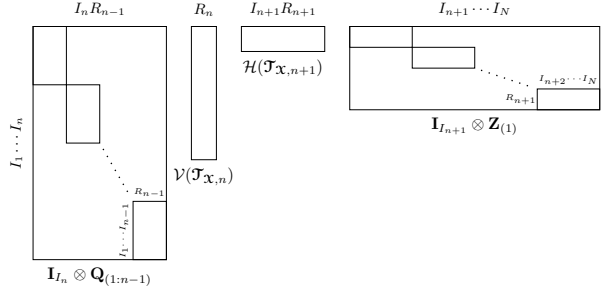

In order to truncate a particular rank $R_n$ by considering only the $n$th TT core using this idea, the TT format should be both $n$-left and $n$-right \change{orthonormal}.
The unfolding of $\T{X}$ that maps the first $n$ tensor dimensions to rows can be expressed as a product of four matrices:
\begin{equation}
\label{eq:quadprodTT}
\Mz{X}{1:n} = (\M{I}_{I_n} \otimes \Mz{Q}{1:n-1}) \cdot \VOp(\TT{X}{n}) \cdot \HOp(\TT{X}{n+1}) \cdot (\M{I}_{I_{n+1}} \otimes \Mz{Z}{1}),
\end{equation}
where $\T Q$ is $I_1 \times \cdots \times I_{n-1} \times R_{n-1}$ with 
$$\T{Q}(i_1,\dots,i_{n-1},r_{n-1}) = \TTslice{X}{1} \cdot \TTslice{X}{2} \cdots \TT{X}{n-1}(:,i_{n-1},r_{n-1}),$$
and $\T Z$ is $R_{n+1} \times I_{n+2} \times \dots \times I_N$ with
$$\T{Z}(r_{n+1},i_{n+2},\dots,i_N) = \TT{X}{n+2}(r_{n+1},i_{n+2},:) \cdot \TTslice{X}{n+3} \cdots \TTslice{X}{N}.$$
See \Cref{fig:quadprod} for a visualization and \Cref{app:quadprod} for a full derivation of \cref{eq:quadprodTT}.
If $\T{X}$ is $n$-left and $n$-right \change{orthonormal}, then $\Mz{Q}{1:n-1}$ and $\Mz{Z}{1}$ are column and row \change{orthonormal} (and so are their Kronecker products with an identity matrix), respectively, and $\HOp(\TT{X}{n+1})$ is also row \change{orthonormal}.

In order to truncate $R_n$, we view \cref{eq:quadprodTT} as an instance of \cref{eq:quadprodmatrix} where $\VOp(\TT{X}{n})$ plays the role of $\M{B}$ and $\HOp(\TT{X}{n+1})$ plays the role of $\M{C}$  (though $\HOp(\TT{X}{n+1})$ is already \change{orthonormal}ized).
We compute the truncated SVD $\VOp(\TT{X}{n})\approx \M[\hat]{U}\M[\hat]{\Sigma}\M[\hat]{V}^\Tra$, replace $\VOp(\TT{X}{n})$ with $\M[\hat]{U}$, and apply $\M[\hat]{\Sigma}\M[\hat]{V}^\Tra$ to $\HOp(\TT{X}{n+1})$.
In this way, $R_n$ is truncated, $\VOp(\TT{X}{n})$ becomes column \change{orthonormal}, and because $\T{Q}$ and $\T{Z}$ are not modified, $\T{X}$ becomes $(n{+}1)$-left and $(n{+}1)$-right \change{orthonormal} and ready for the truncation of $R_{n+1}$. 

The rounding procedure consists of two sweeps along the modes. 
During the first, the tensor is left or right \change{orthonormal}ized. 
On the second, sweeping in the opposite direction, the TT ranks are reduced sequentially via SVD truncation of the matricized cores.
The rounding accuracy $\varepsilon_0$ can be defined a priori such that the rounded TT tensor is $\varepsilon_0$-close to the original TT tensor.
We note that this method is quasi-optimal in finding the closest TT tensor with prescribed TT ranks to a given TT tensor \cite{OT10}.

\begin{algorithm}
\caption{TT-rounding}
\label{alg:TT-rounding}
\begin{algorithmic}[1]
	\Require{A tensor $\T{Y}$ in TT format, a threshold $\varepsilon_0$}
	\Ensure{A tensor $\T{X}$ in TT format with reduced ranks such that \change{$\|\T{X} - \T{Y}\| \leq \varepsilon_0 \|\T{Y}\|$}}
	\Function{$\T{X} =$ Rounding}{$\T{Y}, \varepsilon_0$}
		\State $\T{X}$ = \Call{Right-\change{Orthonormal}ization}{$\T{Y}$}
		\State{Compute $\|\T{Y}\| = \|\TT{X}{1}\|_F$ and the truncation threshold $\varepsilon = \frac{\|\T{Y}\|}{\sqrt{N-1}} \varepsilon_0$}
		\For{$n=1$ to $N-1$}
			\State{$ [\VOp(\TT{X}{n}), \M{\Sigma}, \M{V}] = \text{SVD}(\VOp(\TT{X}{n}),\varepsilon)$ \Comment{$\varepsilon$-truncated SVD factorization}}
			\State{$\HOp(\TT{X}{n+1}) = \M{\Sigma} \M{V}^\Tra \HOp(\TT{X}{n+1})$ \Comment{$\TT{X}{n+1}=\TT{X}{n+1} \times_1 (\M{\Sigma}\M{V}^\Tra)$}}
		\EndFor
	\EndFunction
\end{algorithmic}
\end{algorithm}

\subsection{\change{Parallel Cost Model}}
\label{sec:par-model}

\change{%
To analyze our parallel algorithms, we use the MPI-based model that tracks floating point operations (flops) as well as the amount of data and number of messages communicated along the critical path \cite{BCDH+14,CH+07,TRG05}.
In this model, communication is performed via point-to-point messages, and the time is estimated as the sum of time spent in computation and communication along the critical path.
In this way, processors can perform independent computations simultaneously and disjoint pairs of processors can communicate messages simultaneously.
Each flop is assumed to cost $\gamma$ units of time, and message of $n$ words is assumed to cost $\alpha+\beta n$ units, where $\alpha$ is referred to as the per-message latency cost and $\beta$ is the per-word bandwidth cost.
Accumulating costs along the critical path ensures that computation and communication that depend on one another occur in sequence.
The parallel time cost is thus estimated as $\gamma \cdot \text{\# flops} + \beta \cdot \text{\# words} + \alpha \cdot \text{\# messages}$.
Overlapping computation and communication is a useful optimization (and our implementation does so when possible), but the model ignores this possibility as it affects the overall running time by at most a constant.
Algorithms for collective communications among groups of processors, such as \textsc{AllReduce}, have been optimized for this model (and within MPI implementations), and we use the previously established costs of collectives \cite{CH+07,TRG05} in our analysis.
}

\section{Parallel Algorithms for Tensor Train}
\label{sec:parallel-tt}

In this section we detail the parallel algorithms for manipulating TT tensors that are distributed over multiple processors' memories.
We describe our proposed data distribution of the core tensors in \Cref{sec:distribution}, which is designed for efficient \change{orthonormal}ization and truncation of TT tensors.
In \Cref{sec:basicops} we show how to perform basic operations on TT tensors in this distribution such as addition, elementwise multiplication, and applying certain linear operators.
Our proposed parallel \change{orthonormal}ization and truncation routines are presented in \Cref{sec:par_orth,sec:par_round}, respectively.
Both of those routines rely on an existing communication-efficient parallel QR decomposition algorithm called Tall-Skinny QR (TSQR) \cite{DGHL12}, which is given for completeness in \Cref{sec:TSQR}.
A summary of the costs of the parallel algorithms is presented in \Cref{tab:costs}.

\begin{table}
\begin{tabular}{|c|ccc|}
\hline
\textbf{TT Algorithm} & \textbf{Computation} & \textbf{Comm.~Data} & \textbf{Comm.~Msgs} \\ 
\hline
Summation & --- & --- & --- \\
Hadamard & $\frac{NIR^4}{P}$ & --- & --- \\
Inner Product & $4\frac{NIR^3}{P}$ & $O(NR^2)$ & $O(N\log P)$ \\ 
Norm & $2\frac{NIR^3}{P}$ & $O(NR^2)$ & $O(N\log P)$ \\
\change{Orthonormal}ization & $5\frac{NIR^3}{P} + O(NR^3\log P)$ & $O(NR^2\log P)$ & $O(N\log P)$ \\
Rounding & $7\frac{NIR^3}{P} + O(NR^3\log P)$ & $O(NR^2\log P)$ & $O(N\log P)$ \\
\hline
\end{tabular}
\caption{Summary of computation and communication costs of parallel TT operations using $P$ processors, assuming inputs are $N$-way tensors with identical dimensions $I_n=I$ and ranks $R_n=R$. The computation cost of rounding assumes the original ranks are reduced in half; the constant can range from 3 to 13 depending on the reduced ranks.}
\label{tab:costs}
\end{table}

\subsection{Data Distribution and Layout}
\label{sec:distribution}

We are interested in the parallelization of TT operations with a large number of modes and where one or multiple mode sizes are very large comparing to the TT ranks.
This type of configuration arises in many applications such as parameter dependent PDEs \cite{KreT10}, stochastic PDEs \cite{KreKNT15}, and molecular simulations \cite{SavDWK14}.
In case there exist TT cores with relatively small mode sizes, those can be stored redundantly on each processor.
We note that our implementation can deal with both cases.
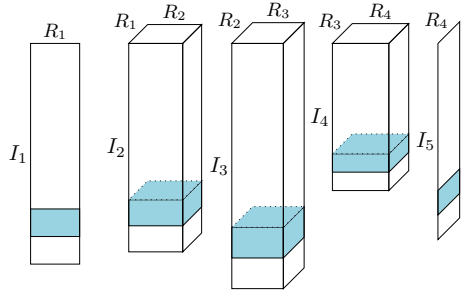
\begin{figure}
	\begin{center}
  \begin{adjustbox}{height=0.3\linewidth}

\pgfmathsetmacro{\np}{8} 
\pgfmathsetmacro{\pid}{6} 

\begin{tikzpicture}[scale=.8,textnode/.style={scale=.5}]
\pgfmathsetmacro{\ione}{(\pid/\np)*\none}
\pgfmathsetmacro{\vecwidth}{(1/\np)*\none}
\draw[fill=mpiblue!50] (0,-\ione,0) rectangle (-\rtwo,-\ione-\vecwidth,0);
\begin{scope}[shift={(\rthree+\gap,0,0)}]
	\pgfmathsetmacro{\itwo}{(\pid/\np)*\ntwo}
	\pgfmathsetmacro{\vecwidth}{(1/\np)*\ntwo}
	\draw[fill=mpiblue!50] (0,-\itwo,0) rectangle (-\rthree,-\itwo-\vecwidth,0);
	\draw[fill=mpiblue!50] (0,-\itwo,0) -- ++(0,0,-\rtwo) -- ++(0,-\vecwidth,0) -- ++(0,0,\rtwo) -- cycle;
	\draw[fill=mpiblue!50,dotted] (0,-\itwo,0) -- ++(-\rthree,0,0) -- ++(0,0,-\rtwo) -- ++(\rthree,0,0) -- cycle;
\end{scope}
\begin{scope}[shift={(\rthree+\rfour+2*\gap,0,0)}]
	\pgfmathsetmacro{\ithree}{(\pid/\np)*\nthree}
	\pgfmathsetmacro{\vecwidth}{(1/\np)*\nthree}
	\draw[fill=mpiblue!50] (0,-\ithree,0) rectangle (-\rfour,-\ithree-\vecwidth,0);
	\draw[fill=mpiblue!50] (0,-\ithree,0) -- ++(0,0,-\rthree) -- ++(0,-\vecwidth,0) -- ++(0,0,\rthree) -- cycle;
	\draw[fill=mpiblue!50,dotted] (0,-\ithree,0) -- ++(-\rfour,0,0) -- ++(0,0,-\rthree) -- ++(\rfour,0,0) -- cycle;
\end{scope}
\begin{scope}[shift={(\rthree+\rfour+\rfive+3*\gap,0,0)}]
	\pgfmathsetmacro{\ifour}{(\pid/\np)*\nfour}
	\pgfmathsetmacro{\vecwidth}{(1/\np)*\nfour}
	\draw[fill=mpiblue!50] (0,-\ifour,0) rectangle (-\rfive,-\ifour-\vecwidth,0);
	\draw[fill=mpiblue!50] (0,-\ifour,0) -- ++(0,0,-\rfour) -- ++(0,-\vecwidth,0) -- ++(0,0,\rfour) -- cycle;
	\draw[fill=mpiblue!50,dotted] (0,-\ifour,0) -- ++(-\rfive,0,0) -- ++(0,0,-\rfour) -- ++(\rfive,0,0) -- cycle;
\end{scope}
\begin{scope}[shift={(\rthree+\rfour+\rfive+4*\gap,0,0)}]
	\pgfmathsetmacro{\ifive}{(\pid/\np)*\nfive}
	\pgfmathsetmacro{\vecwidth}{(1/\np)*\nfive}
	\draw[fill=mpiblue!50] (0,-\ifive,0) -- ++(0,0,-\rfive) -- ++(0,-\vecwidth,0) -- ++(0,0,\rfive) -- cycle;
\end{scope}
\drawTT 
\end{tikzpicture}
    \end{adjustbox}
	\end{center}
  \caption{1D distribution of a TT tensor across $P$ processors with data owned by a particular processor highlighted in blue.}
  \label{fig:distributed_core}
\end{figure}

Algorithms for \change{orthonormal}ization and rounding of TT tensors are sequential with respect to the mode; often computation can occur on only one mode at a time.
In order to utilize all processors and maintain load balancing in a parallel environment, we choose to distribute each TT core over all processors, so that each processor owns a subtensor of each TT core.
To ensure the computations on each core can be done in a communication-efficient way, we choose a 1D distribution for each core, where the mode corresponding to the original tensor is divided across processors.
This corresponds to a Cartesian distribution of each $R_{n-1}\times I_n\times R_n$ core over a $1 \times P \times 1$ processor grid, or equivalently, a block row distribution of $\VOp(\TT{X}{n})$ or a block column distribution of $\HOp(\TT{X}{n})$, for $n=1,\ldots,N$; see~\Cref{fig:distributed_core}.
In this manner, each processors owns $N$ local subtensors with dimensions $\{R_{n-1}\times (I_n/P) \times R_n\}$.
The notation $\TT{X}{n}^{(p)}$ denotes the local subtensor of the $n$th core owned by processor $p$.

This distribution allows performing basic operations, such as addition and elementwise multiplication, on the TT representation locally, see~\Cref{sec:basicops}.
\change{Furthermore, the bottleneck computations within orthonormalization and rounding are orthonormalization of vertical and horizontal unfoldings of TT cores. 
For communication optimality of these operations, the TSQR algorithm (see~\Cref{sec:TSQR}) requires that both of these unfoldings are in 1D matrix distribution, which in turn requires that the TT core be distributed over a $1 \times P \times 1$ processor grid.}
The distribution of a TT core in this way can also be seen as a generalization of the distribution of a vector in parallel iterative linear solvers \cite{Ald18,Jol14}.
Indeed, if $\M{A}$ is an $I_n \times I_n$ sparse matrix distributed across processors as block row panels, the computation of $\M{A} \TT{X}{n}(k, :, l)$ can be done by using \change{standard parallel} sparse-matrix-vector multiplication routines.
\change{We note that a drawback of this distribution is that the available parallelism in each TT core computation is limited to the size of the tensor dimension.
If the TT ranks are much larger than the tensor dimension, then alternative distributions, redistributions, and parallelizations should be considered.}


Tensor entries are linearized in memory.
Each local core tensor $\TT{X}{n}^{(p)}$ is $R_{n-1}\times (I_n/P) \times R_n$, and we store it in the ``vec-oriented'' or ``natural descending'' order \cite{BalKK20,RVL12} in memory.
For 3-way tensors, this means that mode-$1$ fibers (of length $R_{n-1}$) are contiguous in memory, as this corresponds to the mode-$1$ modal unfolding.
Additionally, the mode-$3$ slices (of size $R_{n-1}\times (I_n/P)$) are also contiguous in memory and internally linearized in column-major order, as this corresponds to the more general $(1{:}2)$ unfolding \cite{PTC13a,RVL12}.
In particular, these facts imply that both the vertical and horizontal unfoldings are column major in memory. 

BLAS and LAPACK routines require either row- or column-major ordering (unit stride for one dimension and constant stride for the other), but this property of the vertical and horizontal unfoldings means that we can operate on them without any physical permutation of the tensor data. 
For example, we can perform operations such as QR factorization of $\VOp(\TT{X}{n})$ and $\VOp(\TT{X}{n}) \M{R}$, where $\M{R} \in \R^{R_n \times R_n}$, with a single LAPACK or BLAS call.

This choice of ordering comes at the expense of less convenient access to the mode-2 modal unfolding (of dimension $(I_n/P) \times R_{n-1}R_n$), which is neither row or column major in memory.
This unfolding can be visualized in memory as a concatenation of $R_n$ contiguous submatrices, each of dimension $(I_n/P) \times R_{n-1}$ and each stored in row-major order \cite{BalKK20}.
In order to perform the mode-2 multiplication (tensor times matrix operation), as is necessary in the application of a spatial operator on the core, we must make a sequence of calls to the matrix-matrix multiplication BLAS subroutine.
That is, we make $R_n$ calls for multiplications of the same $I_n \times I_n$ matrix with different $I_n \times R_{n-1}$ matrices.

\subsection{Basic Operations}
\label{sec:basicops}

\subsubsection{Summation}
\label{sec:sum}

To sum two tensors $\T{X}$ and $\T{Y}$, we can write \cite{Ose11}:
\begin{align*}
\Tentry{Z} ={}& \Tentry{X} + \Tentry{Y} \\
	={}& \TTslice{X}{1} \cdots \TTslice{X}{N} + \TTslice{Y}{1} \cdots \TTslice{Y}{N} \\
\begin{split}          
	={}& \begin{pmatrix} \TTslice{X}{1} & \TTslice{Y}{1} \end{pmatrix} \begin{pmatrix} \TTslice{X}{2} & \\ & \TTslice{Y}{2} \end{pmatrix} \\
	&{} \cdots \begin{pmatrix} \TTslice{X}{N-1} & \\ & \TTslice{Y}{N-1} \end{pmatrix}
	\begin{pmatrix} \TTslice{X}{N} \\ \TTslice{Y}{N} \end{pmatrix}.
\end{split}
\end{align*}
Thus, the TT representation of $\T{Z} = \T{X} + \T{Y}$ is given by the following slice-wise formula: 
$$\TTslice{Z}{n} = \begin{pmatrix} \TTslice{X}{n} & \\ & \TTslice{Y}{n} \end{pmatrix}$$
for $2\leq n\leq N-1$, and $1\leq i_n\leq I_n$.
We also have $\TT{Z}{1} =  \begin{pmatrix} \TT{X}{1} & \TT{Y}{1} \end{pmatrix}$ and $\TT{Z}{N} = \begin{pmatrix} \TT{X}{N} \\ \TT{Y}{N} \end{pmatrix}$.
\change{Note that the TT ranks of this representation of $\T{Z}$ are the sums of the TT ranks of $\T{X}$ and $\T{Y}$.}

Given the 1D data distribution of each core described in \Cref{sec:distribution}, the summation operation can be performed locally with no interprocessor communication.
That is, because $\T X$, $\T Y$, and $\T Z$ have identical dimensions, they will have identical distributions, and each slice of a core tensor of $\T Z$ will be owned by the processor that owns the corresponding slices of cores of $\T X$ and $\T Y$.

\subsubsection{Hadamard Product}

To compute the Hadamard (elementwise) product of two tensors $\T X$ and $\T Y$, we can write \cite{Ose11}:
\begin{align*}
\Tentry{Z} ={}& \Tentry{X} \cdot \Tentry{Y} \\
	={}& \left(\TTslice{X}{1} \cdots \TTslice{X}{N} \right)\cdot \left(\TTslice{Y}{1} \cdots \TTslice{Y}{N}\right) \\
	={}& \left(\TTslice{X}{1} \cdots \TTslice{X}{N} \right)\otimes \left(\TTslice{Y}{1} \cdots \TTslice{Y}{N}\right) \\
	={}& \left(\TTslice{X}{1} \otimes \TTslice{Y}{1}\right) \cdots \left(\TTslice{X}{N} \otimes \TTslice{Y}{N}\right). 
\end{align*}
Thus, the TT representation of $\T{Z} = \T{X} \Hada \T{Y}$ is given by the following slice-wise formula: 
$\TTslice{Z}{n} = \TTslice{X}{n} \otimes \TTslice{Y}{n}$ for $1\leq n\leq N$ and $1\leq i_n\leq I_n$.
Here, the \change{TT ranks of the representation of} $\T{Z}$ are the products of the TT ranks of $\T{X}$ and $\T{Y}$.

Again, given the 1D data distribution of each core and the fact that each core is computed slice-wise, the Hadamard product can be performed locally with no interprocessor communication.
We note that because of the extra expense of the Hadamard product (due to computing explicit Kronecker products of slices), it is likely advantageous to maintain Hadamard products in implicit form for later operations such as rounding.
\change{While we do not pursue this approach further in this work,} the combination of Hadamard products and recompression has been shown to be effective for Tucker tensors, \change{but it requires randomization in the truncation operations} \cite{KP17}.

\subsubsection{Inner Product}
\label{sec:inner}

To compute the inner product of two tensors $\T X$ and $\T Y$, using similar identities as for the Hadamard product, we can write \cite{Ose11}:
\begin{align*}
\left \langle \T X, \T Y \right \rangle ={}& \sum_{i_1,\ldots,i_N} \Tentry{X} \cdot \Tentry{Y} \\
	={}& \sum_{i_1,\ldots,i_N}\left(\TTslice{X}{1} \otimes \TTslice{Y}{1}\right) \cdots \left(\TTslice{X}{N} \otimes \TTslice{Y}{N}\right) \\
\begin{split}
	={}& \sum_{i_1} \left(\TTslice{X}{1} \otimes \TTslice{Y}{1}\right) \sum_{i_2} \left( \TTslice{X}{2} \otimes \TTslice{Y}{2}\right) \\ 
	&{} \cdots \sum_{i_N} \left(\TTslice{X}{N} \otimes \TTslice{Y}{N}\right). 
\end{split}
\end{align*}
This expression can be evaluated efficiently by a sequence of structured matrix-vector products that avoid forming Kronecker products of matrices, and these matrix-vector products are cast as matrix-matrix multiplications.

To see how, we assume that the TT ranks of $\T X$ and $\T Y$ are $\{R^{\T X}_n\}$ and $\{R^{\T Y}_n\}$, respectively.
First, we explicitly construct the row vector 
$$\V{w}_1 = \sum_{i_1}\TTslice{X}{1} \otimes \TTslice{Y}{1},$$
which has dimension $R^{\T X}_1\cdot R^{\T Y}_1$.
Note that $\V{w}_1$ is the vectorization of the matrix $\VOp(\TT{Y}{1})^\top \VOp(\TT{X}{1})$.
Then we distribute $\V{w}_1$ to all terms within the next summation to compute $\V{w}_2$ using
$$\V{w}_2 = \sum_{i_2} \V{w}_1 \left( \TTslice{X}{2} \otimes \TTslice{Y}{2} \right),$$
with each term in the summation evaluated via $\texttt{vec}\left(\TTslice{Y}{2}^\Tra \M{W}_1 \TTslice{X}{2}\right)$, where $\M{W}_1$ is a reshaping of the vector $\V{W}_1$ into a $R^{\T Y}_1\times R^{\T X}_1$ matrix and \texttt{vec} is a row-wise vectorization operator.
We note that $\TTslice{X}{2}$ is $R^{\T X}_1\times R^{\T X}_2$, and $\TTslice{Y}{2}$ is $R^{\T Y}_1\times R^{\T Y}_2$, and $\V{w}_2$ therefore has dimension $R^{\T X}_2\cdot R^{\T Y}_2$.
This process is repeated with
\begin{equation}
\label{eq:innerprodsum}
\M{W}_n = \sum_{i_n} \TTslice{Y}{n}^\Tra \M{W}_{n-1} \TTslice{X}{n},
\end{equation}
until the last core, when we compute the inner product as 
  $\left \langle \T X, \T Y \right \rangle = \sum_{i_N} \TTslice{Y}{N}^\Tra \; \M{W}_{N-1} \; \TTslice{X}{N}$,
where $\M{W}_{N-1}$ is a $R^{\T Y}_{N-1}\times R^{\T X}_{N-1}$ matrix.

If all the tensor dimensions are the same and all TT ranks are the same, i.e., $I=I_1 = \cdots = I_N$ and $R = R^{\T X}_1 = R^{\T Y}_1 = \cdots = R^{\T X}_{N-1} = R^{\T Y}_{N-1}$, the computational complexity is approximately $4NIR^3$.

Evaluating \cref{eq:innerprodsum} directly can exploit the efficiency of dense matrix multiplication, but it requires many calls to the BLAS subroutine.
With some extra temporary memory, we can reduce the number of BLAS calls to 2, performing the same overall number of flops.
Let $\T Z$ be defined such that $\HOp(\TT{Z}{n}) = \M{W}_{n-1} \HOp(\TT{X}{n})$, or the mode-1 multiplication between the core and the matrix, for $n = 1, \ldots, N$ (with $\M{W}_0 = 1$). 
Then, we have $\M{W}_n$ as a contraction of modes 1 and 2 between cores of $\T{Y}$ and $\T{Z}$, or 
$$\M{W}_n = \VOp(\TT{Y}{n})^\top \VOp(\TT{Z}{n}), \quad \text{ for } \ n = 1, \ldots, N.$$
Each of these two multiplications requires a single BLAS call because horizontal and vertical unfoldings are column major in memory.
We note the final contraction in mode $N$ is a dot product instead of a matrix multiplication.

When the input TT tensors are distributed across processors as described in \Cref{sec:distribution}, we can compute the inner product using this technique.
Each term in the summation of \cref{eq:innerprodsum}, which involves corresponding slices of the input tensors, is evaluated by a single processor as long as the matrix $\M{W}_n$ is available on each processor.
Thus, the computation can be load balanced across processors as long as the distribution is load balanced, and each processor can apply the optimization to reduce BLAS calls independently.
We perform an \textsc{AllReduce} collective operation to compute the summation for each mode.
With constant tensor dimensions and TT ranks, the computational cost is approximately $4NIR^3/P$ and the communication cost is $\beta\cdot O(NR^2) + \alpha\cdot O(N\log P)$.

\subsubsection{Norms}
\label{sec:norms}

To compute the norm of a tensor in TT format, we consider two approaches.
The first approach is to use the inner product algorithm described in \Cref{sec:inner} and the identity $\|\T X\|^2 = \langle \T X, \T X \rangle$.
We note that in this case, the matrices $\{\M{W}_n\}$ are symmetric and positive semi-definite, see~\cref{eq:innerprodsum}, and the structured matrix-vector products can exploit this property to save roughly half the computation.
Since $\M{W}_n$ is SPSD, it admits a triangular factorization given by pivoted Cholesky (or LDL): $\M{W}_n = \M{P}_n\M{L}_n \M{L}_n^\top \M{P}_n^\top$. 
Thus, the matrix $\M{W}_{n}$ is computed as $\M{W}_{n} = \VOp(\TT{Z}{n})^\top \VOp(\TT{Z}{n})$, where $\HOp(\TT{Z}{n}) = \M{L}_{n-1}^\top (\M{P}_{n-1}^\top \HOp(\TT{X}{n}))$.
The triangular multiplication to compute the $n$th core of $\T{Z}$ and the symmetric multiplication to compute $\M{W}_n$ each require half the flops of a normal matrix multiplication, so the overall computational complexity of this approach is $2NIR^3$.
It is parallelized similarly to the general inner product.

The second approach is to first right- or left-\change{orthonormal}ize the tensor using \Cref{alg:TT-rorthonormalization}, and then the norm of the tensor is given by $\|\TT{X}{1}\|_F$ or $\|\TT{X}{N}\|_F$ as shown in \Cref{sec:orth}.
\change{This approach can be more accurate than the first one when computing small norms, as the first approach can suffer from cancellation error.}
When the TT tensor is distributed, the \change{orthonormal}ization procedure is more complicated than computing inner products; we describe the parallel algorithm in \Cref{sec:par_orth}.

\subsubsection{Matrix-Vector Multiplication}

In order to build Krylov-like iterative methods to solve linear systems with solutions in TT-format, we must also be able to apply a matrix operator to a vector in TT-format.
We will consider a restricted set of matrix operators: sums of Kronecker products of \change{sparse} matrices \cite{BM02,KreKNT15,KreT10,Tyr03}.

Each term in the sum can be seen as a generalization of a rank-one tensor to the operator case.
We use the notation 
$$\M{A} = \M{A}_1 \otimes \cdots \otimes \M{A}_N$$
to denote a single Kronecker product of matrices, where the dimensions of $\M{A}_n$ are $I_n\times I_n$, conforming to the dimensions of $\T X$ in TT-format.
In this case, we can compute the matrix-vector multiplication $\text{vec}(\T{Y}) = \M{A} \cdot \text{vec}(\T{X})$, where 
\begin{align*}
\Tentry{Y} ={}& \sum_{j_1,\dots,j_N} \M{A}_1(i_1,j_1) \cdots \M{A}_N(i_N,j_N) \cdot \T{X}(j_1,\dots,j_N) \\
={}& \sum_{j_1,\dots,j_N} \M{A}_1(i_1,j_1) \cdots \M{A}_N(i_N,j_N) \cdot \TT{X}{1}(j_1,:) \cdots \TT{X}{N}(:,j_N) \\
={}& \sum_{j_1} \M{A}_1(i_1,j_1)\TT{X}{1}(j_1,:) \cdots \sum_{j_N} \M{A}_N(i_N,j_N) \TT{X}{N}(:,j_N) \\
={}& \TTslice{Y}{1} \cdots \TTslice{Y}{N}
\end{align*}
with $\TT{Y}{1}=\M{A}_1\TT{X}{1}$, $\TT{Y}{n}=\TT{X}{n} \times_2 \M{A}_n$ for $1<n<N$, and $\TT{Y}{N}=\TT{X}{N}\M{A}_N^\Tra$. 
Here the notation $\times_2$ refers to the mode-2 tensor-matrix product, defined so that 
$$\TT{Y}{n}(r_{n-1},:,r_n) = \M{A}_n\TT{X}{n}(r_{n-1},:,r_n)$$
for $1<n<N$, $1\leq r_{n-1} \leq R_{n-1}$, and $1\leq r_n\leq R_n$.

Thus, applying a Kronecker product of matrices to a vector in TT-format maintains the TT-format with the same ranks, and operations on cores can be performed independently.
In order to apply an operator that is a sum of multiple Kronecker products of matrices, we can apply each term separately and use the summation procedure described in \Cref{sec:sum} along with TT-rounding to control rank growth.
We note that it  is possible to apply more general forms of tensorized operators to vectors in TT-format \cite{Ose11}, but we do not consider them here.

When the vector in TT-format is distributed as described in \Cref{sec:distribution}, we must perform the mode-2 tensor-matrix product using a parallel algorithm.
We can view the mode-2 tensor-matrix product as applying the matrix to the mode-2 unfolding of the tensor core $\TT{X}{n}$ (often denoted with subscript $(2)$ \cite{KolB09}), which has dimensions $I_n\times R_{n-1}R_n$.
We observe that the parallel distribution of the mode-2 unfolding of $\TT{X}{n}$ is 1D row-distributed: each processor owns a subset of the rows of the matrix (corresponding to slices of the core tensor).
Thus, the application of $\M{A}_n$ to this unfolding has the same algorithmic structure as the sparse-matrix-times-multiple-vectors operation (SpMM) where all vectors have the same parallel distribution.
Assuming the matrix $\M{A}_n$ is sparse and also row-distributed, as is common in libraries such as PETSc \cite{PETSc} and Trilinos \cite{HB+05}, the parallel algorithm involves communication of input tensor core slices among processors, where the communication pattern is determined by $\M{A}_n$ and its distribution. 
We do not explore experimental results for such matrix-vector multiplications in this paper, as the performance depends heavily on the application and sparsity structure of the operator matrices.

\subsection{TSQR}
\label{sec:TSQR}

\change{As is evident in \cref{alg:TT-rorthonormalization,alg:TT-rounding}, the QR factorization of tall-skinny matrices is a key subroutine in TT rounding.}
To compute the QR factorizations within the \change{TT orthonormalization and TT rounding procedures} in parallel, we use the Tall-Skinny QR algorithm \cite{DGHL12}, which is designed (and communication efficient) for matrices with many more rows than columns.
For completeness, we present the TSQR subroutine as \cref{alg:BFTSQR}, which corresponds to \cite[Alg. 7]{BDGJ+15}, and the TSQR-Apply-Q subroutine as \cref{alg:BFappQ}.
\change{While TSQR is strictly a matrix algorithm, it is fundamental to the TT algorithms and analysis of \Cref{sec:par_orth,sec:par_round}, so we present it separately in this subsection.}
The subroutines assume a power-of-two number of processors to simplify the pseudocode; see \cref{app:tsqr} for the generalizations to any number of processors.

For a tall-skinny matrix that is 1D row distributed over processors (as is the case for the vertical unfolding and the transpose of the horizontal unfolding), the parallel Householder QR algorithm requires synchronizations for each column of the matrix (to compute and apply each Householder vector).
\change{Furthermore, the local computation of Householder QR is nearly always memory-bandwidth bound in the form of BLAS-2 subroutines (matrix-vector operations).}
The idea of the TSQR algorithm is that the entire factorization can be computed using a single reduction across processors\change{, and each local computation becomes a smaller QR factorization}.
\change{That is, while parallel Householder QR has latency cost of $O(b)$ for a matrix with $b$ columns, TSQR has latency cost $O(\log P)$ (see \Cref{sec:TSQR-fact}).
The superior performance of TSQR over Householder QR has been demonstrated on both distributed-memory and shared-memory platforms \cite{ABDK11,BDGJ+15,CG11,MHDY09}.}

The price of TSQR is that the implicit representation of the \change{orthonormal} factor is more complicated than a single set of Householder vectors, and that the representation depends on the structure of the reduction tree.
We can maintain and apply the \change{orthonormal} factor in this implicit form as long as the parallel algorithm for applying it uses a consistent tree structure.
We note that we employ the ``butterfly'' variant of TSQR, which corresponds to an \change{\textsc{AllReduce}}-like collective operation such that at the end of the algorithm the triangular factor $\M{R}$ is owned by all processors redundantly.
\change{At each of the $\log P$ steps, each processor determines a different partner processor with which to exchange data.}
Another variant uses a binomial tree, corresponding to a reduce-like collective with the triangular factor owned by a single processor.
\change{In the context of TT, the key advantage of the butterfly over the binomial variant is the reduction in communication when the implicit orthogonal factor is applied to another matrix, as we describe in \Cref{sec:TSQR-apply}.}
We compare performance of these two variants in \Cref{sec:perf_tsqr}.

\subsubsection{Factorization}
\label{sec:TSQR-fact}

TSQR (\cref{alg:BFTSQR}) has two phases: \change{local submatrix orthonormal}ization (\cref{line:BFTSQR-leaf}) and parallel reduction of remaining triangular factors (\cref{line:BFTSQR-tree-start} through \cref{line:BFTSQR-tree-stop}).
The cost of the TSQR is as follows:
\begin{equation}
\label{eq:tsqrcost}
\gamma\cdot\left(2\frac{mb^2}{P}+O(b^3 \log P)\right) +\beta\cdot O(b^2 \log P) + \alpha\cdot O(\log P),
\end{equation}
where $m$ is the number of rows and $b$ is the number of columns \cite{DGHL12}.
The leading order flop cost is the QR of the local $(m/P)\times b$ submatrix (\cref{line:BFTSQR-leaf}), the leaf of the TSQR tree.
The communication costs come from the TSQR tree, which has height $O(\log P)$.

\begin{algorithm}
\caption{Parallel Butterfly TSQR}
\label{alg:BFTSQR}
\begin{algorithmic}[1]
\Require $\M{A}$ is an $m\times b$ matrix 1D-distributed so that proc $p$ owns row block $\M{A}^{(p)}$
\Require Number of procs is power of two; see \Cref{alg:BFTSQRapp} for general case
\Ensure $\M{A}=\M{Q} \M{R}$ with $\M{R}$ owned by all procs and $\M{Q}$ represented by $\{\M{Y}_{\ell}^{(p)}\}$ with redundancy $\M{Y}_{\ell}^{(p)}=\M{Y}_{\ell}^{(q)}$ for $p \equiv q \mod 2^\ell$ and $\ell < \log P$
\Function{$[\{\M{Y}_{\ell}^{(p)}\}, \M{R}] = $ Par-TSQR}{$\M{A}^{(p)}$}
	\State $p = \textsc{MyProcID}()$
	\State $[\M{Y}_{\log P}^{(p)},\M[\bar]{R}_{\log P}^{(p)}] = \text{Local-QR}(\M{A}^{(p)})$ \Comment{Leaf node QR} \label{line:BFTSQR-leaf}
	\For{$\ell= \log P-1$ down to $0$} \label{line:BFTSQR-tree-start}
		\State $j = 2^{\ell+1} \lfloor \frac{p}{2^{\ell+1}}\rfloor + \change{\left((p + 2^{\ell}) \mod 2^{\ell+1}\right)}$ \Comment{Determine partner}
		\State Send $\M[\bar]{R}_{\ell+1}^{(p)}$ to and receive $\M[\bar]{R}_{\ell+1}^{(j)}$ from proc $j$ \Comment{Communication} \label{line:BFTSQR-swap}
		\If{$p <  j$}
			\State $[\M{Y}_{\ell}^{(p)},\M[\bar]{R}_{\ell}^{(p)}] = \text{Local-QR}\left( \begin{bmatrix} \M[\bar]{R}_{\ell+1}^{(p)} \\ \M[\bar]{R}_{\ell+1}^{(j)} \end{bmatrix} \right)$ \Comment{Tree node QR} \label{line:BFTSQR-internal1}
		\Else
			\State $[\M{Y}_{\ell}^{(p)},\M[\bar]{R}_{\ell}^{(p)}] = \text{Local-QR}\left( \begin{bmatrix} \M[\bar]{R}_{\ell+1}^{(j)} \\ \M[\bar]{R}_{\ell+1}^{(p)} \end{bmatrix} \right)$  \Comment{Partner tree node QR} \label{line:BFTSQR-internal2}
		\EndIf
	\EndFor \label{line:BFTSQR-tree-stop}
	\State $\M{R}=\M[\bar]{R}_{0}^{(p)}$
\EndFunction
\end{algorithmic}
\end{algorithm}


\subsubsection{Applying and Forming $Q$}
\label{sec:TSQR-apply}

The structure of the TSQR-Apply-Q algorithm (\cref{alg:BFappQ}) matches that of TSQR, but in reverse order (because the TSQR algorithm corresponds to applying $Q^\Tra$).
Thus, the root of the tree is applied first and the leaves last.
However, by using a butterfly tree the communication cost of the TSQR-Apply-Q algorithm (\cref{alg:BFappQ}) is $0$ if the number of processors is a power of $2$ and $\beta \cdot bc + \alpha$ otherwise (the cost of one message; see \Cref{app:tsqr}).
The cost of TSQR-Apply-Q is then
\begin{equation}
\label{eq:applyqcost}
\gamma\cdot\left(4\frac{mbc}{P}+O(b^2c \log P)\right) +\beta\cdot bc + \alpha,
\end{equation}
where the additional parameter $c$ is the number of columns of $C$.
The leading order flop cost is the application of the local $Q$ matrix at the leaf of the TSQR tree (\cref{line:BFappQ-leaf}).

Using a binomial tree TSQR algorithm requires more communication in the application phase (see \cite[Algorithm 8]{BDGJ+15}, for example).
We also note that if the input matrix $\M C$ is upper triangular, then the leading constant can be reduced from 4 to 2 by exploiting the sparsity structure in this local application (and within the tree because all $\M[\bar]{B}_\ell^{(p)}$ matrices are upper triangular in this case, throughout the algorithm), which matches the computation cost of the factorization.
In particular, when we form $\M Q$, we use this algorithm with $\M C$ as the identity matrix, which is upper triangular.

\begin{algorithm}
\caption{Parallel Application of Implicit $Q$ from Butterfly TSQR}
\label{alg:BFappQ}
\begin{algorithmic}[1]
\Require $\{\M{Y}_{\ell}^{(p)}\}$ represents \change{orthonormal} matrix $\M{Q}$ computed by \cref{alg:BFTSQR}
\Require $\M{C}$ is $b\times c$ and redundantly owned by all processors
\Require Number of procs is power of two; see \Cref{alg:BFappQapp} for general case

\Ensure $\M{B}=\M{Q}\begin{bmatrix} \M{C} \\ \M{0} \end{bmatrix}$ is $m\times c$ and 1D-distributed so that proc $p$ owns row block $\M{B}^{(p)}$
								\Function{$\M{B} = $ Par-TSQR-Apply-Q}{$\{\M{Y}_{\ell}^{(p)}\},\M{C}$}
\State $p = \textsc{MyProcID}()$
\State $\bar{\M{B}}_0^{(p)} = \M{C}$
\For{$\ell= 0$ to $ \log P -1$}
  \State $j = 2^{\ell+1} \lfloor \frac{p}{2^{\ell+1}}\rfloor + \change{\left((p + 2^{\ell}) \mod 2^{\ell+1}\right)}$ \Comment{Determine partner}
	\If{$p < j$}
		\State $\begin{bmatrix} \bar{\M{B}}_{\ell+1}^{(p)} \\ \M[\bar]{B}_{\ell+1}^{(j)} \end{bmatrix} = \textsc{Loc-Apply-Q}\left(\begin{bmatrix} \M{I}_b \\ \M{Y}_{\ell}^{(p)}\end{bmatrix}, \begin{bmatrix} \M[\bar]{B}_\ell^{(p)} \\ \M{0} \end{bmatrix} \right)$ \Comment{Tree node apply}
				\label{line:BFappQ-internal}
	\Else
		\State $\begin{bmatrix} \M[\bar]{B}_{\ell+1}^{(j)} \\ \M[\bar]{B}_{\ell+1}^{(p)} \end{bmatrix} = \textsc{Loc-Apply-Q}\left(\begin{bmatrix} \M{I}_b \\ \M{Y}_{\ell}^{(p)}\end{bmatrix}, \begin{bmatrix} \M[\bar]{B}_\ell^{(p)} \\ \M{0} \end{bmatrix} \right)$ \Comment{Part.~tree node apply}
	\EndIf
\EndFor
	\State $\M{B}^{(p)} = \textsc{Loc-Apply-Q}\left(\M{Y}_{\log P}^{(p)}, \begin{bmatrix} \M[\bar]{B}_{\log P}^{(p)} \\ \M{0} \end{bmatrix} \right)$ \Comment{Leaf node apply} \label{line:BFappQ-leaf}
\EndFunction
\end{algorithmic}
\end{algorithm}


\subsection{TT \change{Orthonormal}ization}
\label{sec:par_orth}

\change{Given the parallel TSQR algorithm of \Cref{sec:TSQR}, we now present a parallel algorithm for TT Orthonormalization.}
\Cref{alg:par_tt_orthonormalization} shows right \change{orthonormal}ization and is a parallelization of \Cref{alg:TT-rorthonormalization}.
The approach for left \change{orthonormal}ization is analogous.
The algorithm is performed via a sequential sweep over the cores, where at each iteration, an LQ factorization row-\change{orthonormal}izes the horizontal unfolding of a core and the triangular factor is applied to its left neighbor core.
The 1D parallel distribution of each core implies that the transpose of the horizontal unfolding is 1D row distributed, fitting the requirements of the TSQR algorithm.
Note that we perform a QR factorization of the transpose of the horizontal unfolding, which corresponds to an LQ factorization of the unfolding itself.

\Cref{fig:orth} depicts the operations within a single iteration of the sweep.
At iteration $n$, TSQR is applied to the $n$th core in \cref{line:R2LQR} (\cref{fig:orth:qr}) and then the \change{orthonormal} factor is formed explicitly in \cref{line:R2LappQ} (\cref{fig:orth:formq}).
The notation $\{\M{Y}_{\ell,n}^{(p)}\}$ signifies the set of triangular matrices owned by processor $p$ in the implicit representation of the QR factorization of the $n$th core, where $\ell$ refers to the level of the tree and indexes the set.
In the case $P$ is a power of $2$, each processor owns $ \log P $ matrices in its set.
Because the TSQR subroutine ends with all processors owning the triangular factor $\M{R}_n$, each processor can apply it to core $n-1$ in the 3rd mode without further communication via local matrix multiplication in \cref{line:R2LappR} (\cref{fig:orth:appr}).

\Cref{line:R2LQR,line:R2LappQ} have the costs, given by \cref{eq:tsqrcost} and \cref{eq:applyqcost} with $m=I_n R_n$ and $b=c=R_{n-1}$. 
Since the computation to form the explicit $\M Q$ matrix exploits the sparsity structure of the identity matrix the constant $4$ in \cref{eq:applyqcost} is reduced to $2$.
These two lines together cost
$$\gamma\cdot\left(4\frac{I_nR_nR_{n-1}^2}{P}+O(R_{n-1}^3 \log P)\right) +\beta\cdot O(R_{n-1}^2 \log P) + \alpha\cdot O(\log P).$$
\Cref{line:R2LappR} is a local triangular matrix multiplication costing $\gamma \cdot I_{k-1}R_{k-2}R_{k-1}^2/P$.
Assuming $I_k=I$ and $R_k=R$ for $1\leq k\leq N{-}1$, the total cost of TT \change{orthonormal}ization is then
\begin{equation}
\label{eq:R2Lorth}
\gamma \cdot \left(5\frac{NIR^3}{P} + O(NR^3 \log P) \right) + \beta \cdot O(NR^2 \log P) + \alpha \cdot O(N\log P).
\end{equation}

\begin{algorithm}
\caption{Parallel TT-Right-\change{Orthonormal}ization}
\label{alg:par_tt_orthonormalization} 
\begin{algorithmic}[1]
	\Require{$\T{X}$ in TT format with each core 1D-distributed}
	\Ensure{$\T{X}$ is right \change{orthonormal}, in TT format with same distribution}
	\Function{Par-TT-Right-\change{Orthonormal}ization}{$\{\TT{X}{n}^{(p)}\}$}
	\For{$n=N$ down to $2$}
		\State{$[\{\M{Y}_{\ell,n}^{(p)}\},\M{R}_n] = $ \textsc{TSQR}($\HOp(\TT{X}{n}^{(p)})^\Tra$) \Comment{QR factorization}} \label{line:R2LQR}
		\State{$\HOp(\TT{X}{n}^{(p)})^\Tra = $ \textsc{TSQR-Apply-Q}($\{\M{Y}_{\ell,n}^{(p)}\},\M{I}_{R_{n-1}}$) \Comment{Form explicit $\M{Q}$}} \label{line:R2LappQ}
		\State{$\VOp(\TT{X}{n-1}^{(p)}) = \VOp(\TT{X}{n-1}^{(p)}) \cdot {\M{R}_n}^\Tra$ \Comment{Apply $\M{R}$ to previous core}} \label{line:R2LappR}
	\EndFor
	\EndFunction
\end{algorithmic}
\end{algorithm}

\renewcommand{\nthcoredimone}{R_{n{-}2}}
\renewcommand{\nthcoredimthree}{R_{n-1}}
\renewcommand{\nextcoredimone}{R_{n-1}}
\renewcommand{\nextcoredimthree}{R_{n}}

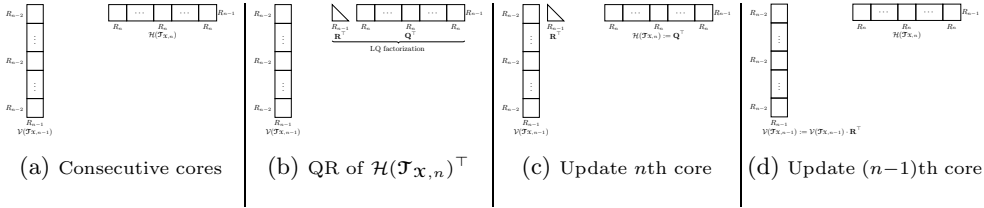
\begin{figure}
  \begin{subfigure}[b]{0.24\linewidth}
    \begin{adjustbox}{width=\linewidth}
      \begin{tikzpicture}[scale=.7,textnode/.style={scale=.5}]
        \begin{scope}[shift={(3*\rtwo,0)}]
          \drawVUnfoldingO
          \node at (0.5*\rtwo,-3.5*\offset) {\small $\VOp(\TT{X}{n-1})$};
        \end{scope}
        \begin{scope}[shift={(8*\rtwo,6*\rfive-1*\rtwo)}]
          \drawHUnfoldingO
          \node at (3*\rthree,-3.5*\offset) {\small $\HOp(\TT{X}{n})$};
        \end{scope}
      \end{tikzpicture}
    \end{adjustbox}
    \captionsub{Consecutive cores}
  \end{subfigure}
  \rulesep
  \begin{subfigure}[b]{0.24\linewidth}
    \begin{adjustbox}{width=\linewidth}
      \begin{tikzpicture}[scale=.7,textnode/.style={scale=.5}]
        \begin{scope}[shift={(3*\rtwo,0)}]
          \drawVUnfoldingO
          \node at (0.5*\rtwo,-3.5*\offset) {\small $\VOp(\TT{X}{n-1})$};
        \end{scope}
        \begin{scope}[shift={(6.5*\rtwo,6*\rfive-1*\rtwo)}]
          \draw (0,0) -- (0, \rtwo) -- (\rtwo, 0) -- cycle;
          \node at (0.5*\rtwo,-1.5*\offset) {\small $\nextcoredimone$};
          \node at (0.5*\rtwo,-3.5*\offset) {\small $\M{R}^\Tra$};
        \end{scope}
        \begin{scope}[shift={(8*\rtwo,6*\rfive-1*\rtwo)}]
          \drawHUnfoldingO
          \node at (3*\rthree,-3.5*\offset) {\small $\M{Q}^\Tra$};
          \draw [thick, decoration = { brace, mirror, raise = 0.5cm}, decorate] (-1.5*\rtwo,-2*\offset) -- (6*\rthree,-2*\offset);
          \node at (2.5*\rtwo, -7*\offset) {\small $\text{LQ factorization}$};
        \end{scope}
      \end{tikzpicture}
    \end{adjustbox}
    \captionsub{QR of $\HOp(\TT{X}{n})^\Tra$} 
    \label{fig:orth:formq}
  \end{subfigure}
  \rulesep
  \begin{subfigure}[b]{0.24\linewidth}
    \begin{adjustbox}{width=\linewidth}
      \begin{tikzpicture}[scale=.7,textnode/.style={scale=.5}]
        \begin{scope}[shift={(3*\rtwo,0)}]
          \drawVUnfoldingO
          \node at (0.5*\rtwo,-3.5*\offset) {\small $\VOp(\TT{X}{n-1})$};
        \end{scope}
        \begin{scope}[shift={(4.5*\rtwo,6*\rfive-1*\rtwo)}]
          \draw (0,0) -- (0, \rtwo) -- (\rtwo, 0) -- cycle;
          \node at (0.5*\rtwo,-1.5*\offset) {\small $\nthcoredimthree$};
          \node at (0.5*\rtwo,-3.5*\offset) {\small $\M{R}^\Tra$};
        \end{scope}
        \begin{scope}[shift={(8*\rtwo,6*\rfive-1*\rtwo)}]
          \drawHUnfoldingO
          \node at (3*\rthree,-3.5*\offset) {\small $\HOp(\TT{X}{n}):=\M{Q}^\Tra$};
        \end{scope}
      \end{tikzpicture}
    \end{adjustbox}
    \captionsub{Update $n$th core}
    \label{fig:orth:qr}
  \end{subfigure}
  \rulesep
  \begin{subfigure}[b]{0.24\linewidth}
    \begin{adjustbox}{width=\linewidth}
      \begin{tikzpicture}[scale=.7,textnode/.style={scale=.5}]
        \begin{scope}[shift={(3*\rtwo,0)}]
          \drawVUnfoldingO
          \node at (0.5*\rtwo+2,-3.5*\offset) {\small $\VOp(\TT{X}{n-1}):=\VOp(\TT{X}{n-1}) \cdot \M{R}^\Tra$};
        \end{scope}
        \begin{scope}[shift={(8*\rtwo,6*\rfive-1*\rtwo)}]
          \drawHUnfoldingO
          \node at (3*\rthree,-3.5*\offset) {\small $\HOp(\TT{X}{n})$};
        \end{scope}
      \end{tikzpicture}
    \end{adjustbox}
    \captionsub{Update $(n{-}1)$th core} 
    \label{fig:orth:appr}
  \end{subfigure}
  \caption{Steps performed in TT right \change{orthonormal}ization}
  \label{fig:orth}
\end{figure}

\subsection{TT Rounding}
\label{sec:par_round}

We present the parallel TT rounding procedure in \Cref{alg:par_tt_rounding}, which is a parallelization of \Cref{alg:TT-rounding}.
The computation consists of two sweeps over the cores, one to \change{orthonormal}ize and one to truncate.
The algorithm shown performs right-\change{orthonormal}ization and then truncates left to right, and the other ordering works analogously.

\begin{algorithm}
  \caption{Parallel TT-Rounding}
  \label{alg:par_tt_rounding} 
  \begin{algorithmic}[1]
          \Require{$\T{X}$ in TT format with each core 1D-distributed over $1{\times} P {\times} 1$ processor grid}
          \Ensure{$\T{Y}$ in TT format with reduced ranks identically distributed across processors}
          \Function{$\{\TT{Y}{n}^{(p)}\} =$ Par-TT-Rounding}{$\{\TT{X}{n}^{(p)}\},\epsilon$}
          \For{$n=N$ down to $2$}
		\State{$[\{\M{Y}_{\ell,n}^{(p)}\},\M{R}_n] = $ \textsc{TSQR}($\HOp(\TT{X}{n}^{(p)})^\Tra$) \Comment{QR factorization}} \label{line:R2LQRi}
		\State{$\VOp(\TT{X}{n-1}^{(p)}) = \VOp(\TT{X}{n-1}^{(p)}) \cdot {\M{R}_n}^\Tra$ \Comment{Apply $\M{R}$ to previous core}} \label{line:R2LappRi}
          \EndFor
          \State Compute $\|\T{X}\|$
          \State $\T{Y}=\T{X}$
          \For{$n=1$ to $N-1$}
		\State{$[\{\M{Y}_{\ell,n}^{(p)}\},\M{R}_n] = $ \textsc{TSQR}($\VOp(\TT{Y}{n}^{(p)})$) \Comment{QR factorization}} \label{line:L2RQRi}
		\State{$[\M[\hat]{U}_R, \hat{\Sigma}, \M[\hat]{V}] = \textsc{tSVD}(\M{R}_n,\frac{\epsilon}{\sqrt{N-1}} \|\T{X}\|)$ \Comment{Redundant truncated SVD of $\M{R}$}} \label{line:L2RSVDi}
		\State{$\VOp(\TT{Y}{n}^{(p)}) = $ \textsc{TSQR-Apply-Q}($\{\M{Y}_{\ell,n}^{(p)}\},\M[\hat]{U}_R$) \Comment{Form explicit $\M[\hat]{U}$}} \label{line:L2RappQi}
		\State{$\HOp(\TT{Y}{n{+}1}^{(p)})^\Tra = $ \textsc{TSQR-Apply-Q}($\{\M{Y}_{\ell,n{+}1}^{(p)}\},\M[\hat]{V} \hat{\Sigma})$ \Comment{Apply $\M[\hat]{\Sigma} \M[\hat]{V}^\Tra$}} \label{line:L2RappVi}
          \EndFor
  \EndFunction
  \end{algorithmic}
\end{algorithm}

\renewcommand{\nthcoredimone}{L_{n{-}1}}
\renewcommand{\nthcoredimthree}{R_{n}}
\renewcommand{\nextcoredimone}{R_{n}}
\renewcommand{\nextcoredimthree}{R_{n{+}1}}

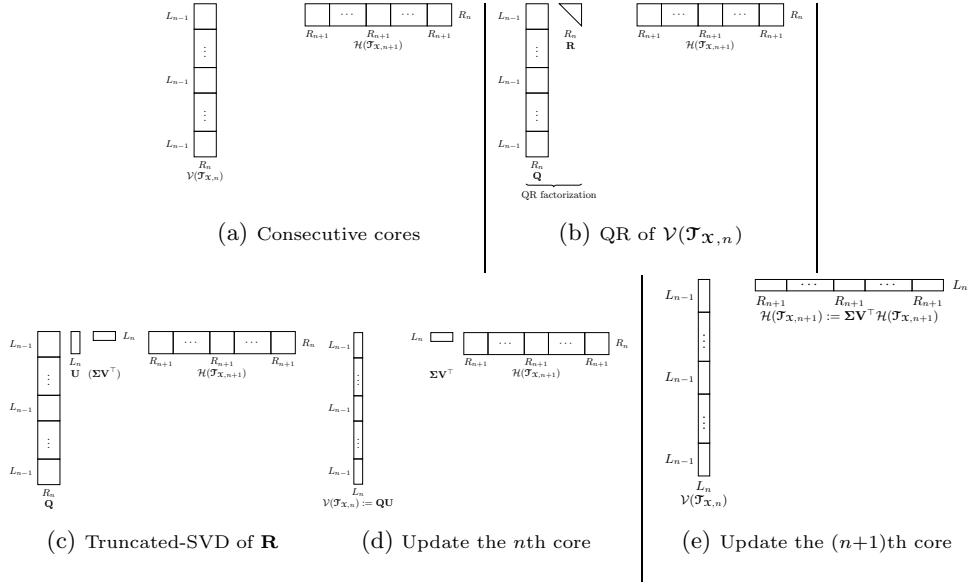
\begin{figure}
  \begin{subfigure}[b]{0.32\linewidth}
    \begin{adjustbox}{width=\linewidth}
      \begin{tikzpicture}[scale=.7,textnode/.style={scale=.5}]
        \begin{scope}[shift={(3*\rtwo,0)}]
          \drawVUnfoldingO
          \node at (0.5*\rtwo,-3.5*\offset) {\small $\VOp(\TT{X}{n})$};
          \node at (1.5*\rtwo, -7*\offset) {\small ${\color{white}\text{QR factorization}}$};
        \end{scope}
        \begin{scope}[shift={(8*\rtwo,6*\rfive-1*\rtwo)}]
          \drawHUnfoldingO
          \node at (3*\rthree,-3.5*\offset) {\small $\HOp(\TT{X}{n+1})$};
        \end{scope}
      \end{tikzpicture}
    \end{adjustbox}
    \captionsub{Consecutive cores}
  \end{subfigure}
  \rulesep
  \begin{subfigure}[b]{0.32\linewidth}
    \begin{adjustbox}{width=\linewidth}
      \begin{tikzpicture}[scale=.7,textnode/.style={scale=.5}]
        \begin{scope}[shift={(3*\rtwo,0)}]
          \drawVUnfoldingO
          \node at (0.5*\rtwo,-3.5*\offset) {\small $\M{Q}$};
          \draw [thick, decoration = { brace, mirror, raise = 0.5cm}, decorate] (0,-2*\offset) -- (2.5*\rtwo,-2*\offset);
          \node at (1.5*\rtwo, -7*\offset) {\small $\text{QR factorization}$};
        \end{scope}
        \begin{scope}[shift={(4.5*\rtwo,6*\rfive-1*\rtwo)}]
          \draw (0,\rtwo) -- (\rtwo, \rtwo) -- (\rtwo, 0) -- cycle;
          \node at (0.5*\rtwo,-1.5*\offset) {\small $\nthcoredimthree$};
          \node at (0.5*\rtwo,-3.5*\offset) {\small $\M{R}$};
        \end{scope}
        \begin{scope}[shift={(8*\rtwo,6*\rfive-1*\rtwo)}]
          \drawHUnfoldingO
          \node at (3*\rthree,-3.5*\offset) {\small $\HOp(\TT{X}{n+1})$};
        \end{scope}
      \end{tikzpicture}
    \end{adjustbox}
    \captionsub{QR of $\VOp(\TT{X}{n})$}
  \end{subfigure}
  \rulesep
  \begin{subfigure}[b]{0.32\linewidth}
    \begin{adjustbox}{width=\linewidth}
      \begin{tikzpicture}[scale=.7,textnode/.style={scale=.5}]
        \begin{scope}[shift={(3*\rtwo,0)}]
          \drawVUnfoldingO
          \node at (0.5*\rtwo,-3.5*\offset) {\small $\M{Q}$};
        \end{scope}
        \begin{scope}[shift={(4.5*\rtwo,6*\rfive-1*\rtwo)}]
          \draw (0,0) rectangle (0.4*\rtwo, \rtwo);
          \node at (0.2*\rtwo,-1.5*\offset) {\small $\nthcoredimthreered$};
          \node at (0.2*\rtwo,-3.5*\offset) {\small $\M{U}$};
        \end{scope}
        \begin{scope}[shift={(5.5*\rtwo,6*\rfive-0.4*\rtwo)}]
          \draw (0,0) rectangle (\rtwo, 0.4*\rtwo);
          \node at (\rtwo + 2.5*\offset,0.2*\rtwo) {\small $\nthcoredimthreered$};
          \node at (0.5*\rtwo,-3.5*\offset-0.6\rtwo) {\small $(\M{\Sigma V}^\top)$};
        \end{scope}
        \begin{scope}[shift={(8*\rtwo,6*\rfive-1*\rtwo)}]
          \drawHUnfoldingO
          \node at (3*\rthree,-3.5*\offset) {\small $\HOp(\TT{X}{n+1})$};
        \end{scope}
      \end{tikzpicture}
    \end{adjustbox}
    \captionsub{Truncated-SVD of $\M{R}$}
  \end{subfigure}
  \begin{subfigure}[b]{0.32\linewidth}
    \begin{adjustbox}{width=\linewidth}
      \begin{tikzpicture}[scale=.7,textnode/.style={scale=.5}]
        \begin{scope}[shift={(3*\rtwo,0)}]
\pgfmathsetmacro{\rtwo}{0.4}
          \drawVUnfoldingOR
          \node at (0.5*\rtwo,-3.5*\offset) {\small $\VOp(\TT{X}{n}) := \M{Q U}$};
        \end{scope}
\pgfmathsetmacro{\rtwo}{1}
        \begin{scope}[shift={(6.5*\rtwo,6*\rfive-0.4*\rtwo)}]
          \draw (0,0) rectangle (\rtwo, 0.4*\rtwo);
          \node at (-2.5*\offset,0.2*\rtwo) {\small $\nthcoredimthreered$};
          \node at (0.5*\rtwo,-0.6*\rtwo-3.5*\offset) {\small $\M{\Sigma V}^\top$};
        \end{scope}
        \begin{scope}[shift={(8*\rtwo,6*\rfive-1*\rtwo)}]
          \drawHUnfoldingO
          \node at (3*\rthree,-3.5*\offset) {\small $\HOp(\TT{X}{n+1})$};
        \end{scope}
      \end{tikzpicture}
    \end{adjustbox}
    \captionsub{Update the $n$th core}
  \end{subfigure}
  \rulesep
  \centering
\pgfmathsetmacro{\rtwo}{0.4}
  \begin{subfigure}[b]{0.32\linewidth}
    \begin{adjustbox}{width=\linewidth}
      \begin{tikzpicture}[scale=.7,textnode/.style={scale=.5}]
        \begin{scope}[shift={(3*\rtwo,0)}]
          \drawVUnfoldingOR
          \node at (0.5*\rtwo,-3.5*\offset) {\small $\VOp(\TT{X}{n})$};
        \end{scope}
        \begin{scope}[shift={(8*\rtwo,6*\rfive-1*\rtwo)}]
          \drawHUnfoldingOR
          \node at (3*\rthree,-3.5*\offset) {\small $\HOp(\TT{X}{n+1}):=\M{\Sigma V}^\top \HOp(\TT{X}{n+1})$};
        \end{scope}
      \end{tikzpicture}
    \end{adjustbox}
    \captionsub{Update the $(n{+}1)$th core} 
  \end{subfigure}
  \caption{Steps performed in iteration of the TT left-to-right truncation}
\end{figure}

\Cref{alg:par_tt_rounding} does not call \Cref{alg:par_tt_orthonormalization} to perform the \change{orthonormal}ization sweep.
This is because \Cref{alg:par_tt_orthonormalization} forms the \change{orthonormal}ized cores explicitly, and \Cref{alg:par_tt_rounding} can leave the \change{orthonormal}ized cores from the first sweep in implicit form to be applied during the second sweep.
 
Iteration $n$ of the right-to-left \change{orthonormal}ization sweep occurs in \Cref{line:R2LQRi,line:R2LappRi}, which matches \Cref{alg:par_tt_orthonormalization} except for the explicit formation of the \change{orthonormal} factor.
Thus, the cost of the \change{orthonormal}ization sweep is 
\begin{equation}
\label{eq:R2LTT}
\gamma \cdot \left(3\frac{NIR^3}{P} + O(NR^3 \log P) \right) + \beta \cdot O(NR^2 \log P) + \alpha \cdot O(N\log P).
\end{equation}

At iteration $n$ of the second loop, \Cref{line:L2RQRi,line:L2RSVDi,line:L2RappQi,line:L2RappVi} implement the left-to-right truncation procedure for the $n$th core in parallel.
\Cref{line:L2RQRi} is a QR factorization and has cost given by \Cref{eq:tsqrcost} with $m=I_nL_{n-1}$ and $b=R_n$, as the number of rows of $\VOp(\TT{Y}{n}^{(p)})$ has been reduced from $I_nR_{n-1}$ to $I_nL_{n-1}$ during iteration $n{-}1$:
$$\gamma\cdot\left(2\frac{I_nL_{n-1}R_n^2}{P}+O(R_n^3 \log P)\right) +\beta\cdot O(R_n^2 \log P) + \alpha\cdot O(\log P).$$
We note that we re-use the notation $\{\M{Y}_{\ell,n}^{(p)}\}$ to store the implicit factorization; while the same variable stored the \change{orthonormal} factor of the $n$th core's horizontal unfolding from the \change{orthonormal}ization sweep, it can be overwritten by this step of the algorithm (the set of matrices will now have different dimensions).
\Cref{line:L2RSVDi} requires $O(R_n^3)$ flops, assuming the full SVD is computed before truncating.
\Cref{line:L2RappQi} implicitly applies an \change{orthonormal} matrix to an $R_n \times L_n$ matrix $\M[\hat]{U}_R$ with cost given by \Cref{eq:applyqcost} with $m=I_nL_{n-1}$, $b=R_n$, and $c=L_n$:
$$\gamma\cdot\left(4\frac{I_nL_{n-1}R_nL_n}{P}+O(R_n^2L_n \log P)\right) +\beta\cdot R_nL_n + \alpha.$$
\Cref{line:L2RappVi} implicitly applies an \change{orthonormal} matrix to an $R_n \times L_n$ matrix $\M[\hat]{V}\M[\hat]{\Sigma}$ with cost given by \Cref{eq:applyqcost} with $m=I_{n+1}R_{n+1}$, $b=R_n$, and $c=L_n$:
$$\gamma\cdot\left(4\frac{I_{n+1}R_{n+1}R_{n}L_n}{P}+O(R_n^2L_n \log P)\right) +\beta\cdot R_nL_n + \alpha.$$

Assuming $I_k=I$, $R_k=R$, and $L_k=L$ for $1\leq k\leq N{-}1$, the total cost of \Cref{alg:par_tt_rounding} is then
\begin{equation}
\label{eq:L2RTT}
\gamma \cdot \left(NIR\frac{3R^2+6RL+4L^2}{P} + O(NR^3\log P) \right) + \beta \cdot O(NR^2 \log P) + \alpha \cdot O(N\log P).
\end{equation}
We note that leaving the \change{orthonormal} factors in implicit form during the \change{orthonormal}ization sweep (as opposed to calling \Cref{alg:par_tt_orthonormalization}) saves up to 40\% of the computation, when the reduced ranks $L_n$ are much smaller than the original ranks $R_n$.
As the rank reduction diminishes, so does the advantage of the implicit optimization.
For example, when ranks are all halved, the reduction in leading order flop cost is 12.5\%.

\section{Numerical Experiments}
\label{sec:numerical_experiments}



\newcommand{\na}{}
\newcommand{\nb}{}
\newcommand{\nc}{}
\newcommand{\nd}{}
\newcommand{\nalabel}{}
\newcommand{\nblabel}{}
\newcommand{\nclabel}{}
\newcommand{\ndlabel}{}
\newcommand{\nelabel}{}
\newcommand{\tra}{}
\newcommand{\trb}{}
\newcommand{\tola}{}
\newcommand{\tolb}{}
\newcommand{\prep}{}
\newcommand{\numiters}{}

\definecolor{color1}{HTML}{2C7BB6}
\definecolor{color2}{HTML}{D7191C}
\definecolor{color3}{HTML}{ABD9E9}
\definecolor{color4}{HTML}{FFFFBF}
\definecolor{color5}{HTML}{FDAE61}
\definecolor{color6}{HTML}{BCFA72}
\definecolor{black}{HTML}{000000}

\newif\iflegend
\legendtrue
\newif\ifylabel
\ylabeltrue

\newcommand{\tsqrplotoptions}{
  ybar stacked,
  reverse legend,
  bar width=10pt,
  \ifylabel
  	ylabel={Time (seconds)}, 
  \fi
  ymin=0,
  ymax=1/\numiters,
  symbolic x coords={TSQR-\tra-\na,TSQR-\trb-\na,,TSQR-\tra-\nb,TSQR-\trb-\nb,,TSQR-\tra-\nc,TSQR-\trb-\nc,,TSQR-\tra-\nd,TSQR-\trb-\nd},
  xticklabels={Binomial,Binomial,Binomial,Binomial,Butterfly,Butterfly,Butterfly,Butterfly},  
  xtick=data,
  x tick label style={rotate=270},
  legend pos={north west},
  legend style={cells={align=left}, nodes={scale=0.8}},
}

\newcommand{\maketsqrplot}{
  \pgfmathsetmacro{\ycoord}{1/\numiters} 
  \begin{axis}[\tsqrplotoptions
    after end axis/.code={ 
      \node(O11) at (axis cs:TSQR-\tra-\na,\ycoord) {};
      \node(O14) at (axis cs:TSQR-\tra-\na,\ycoord) {};
      \node(O1) at ($(O11)!0.5!(O14)$,\ycoord) {};
      \node(N11) at (axis cs:TSQR-\trb-\na,\ycoord) {};
      \node(N14) at (axis cs:TSQR-\trb-\na,\ycoord) {};
      \node(N1) at ($(N11)!0.5!(N14)$,\ycoord) {};
      \node(O21) at (axis cs:TSQR-\tra-\nb,\ycoord) {};
      \node(O24) at (axis cs:TSQR-\tra-\nb,\ycoord) {};
      \node(O2) at ($(O21)!0.5!(O24)$,\ycoord) {};
      \node(N21) at (axis cs:TSQR-\trb-\nb,\ycoord) {};
      \node(N24) at (axis cs:TSQR-\trb-\nb,\ycoord) {};
      \node(N2) at ($(N21)!0.5!(N24)$,\ycoord) {};
      \node(O31) at (axis cs:TSQR-\tra-\nc,\ycoord) {};
      \node(O34) at (axis cs:TSQR-\tra-\nc,\ycoord) {};
      \node(O3) at ($(O31)!0.5!(O34)$,\ycoord) {};
      \node(N31) at (axis cs:TSQR-\trb-\nc,\ycoord) {};
      \node(N34) at (axis cs:TSQR-\trb-\nc,\ycoord) {};
      \node(N3) at ($(N31)!0.5!(N34)$,\ycoord) {};
      \node(O41) at (axis cs:TSQR-\tra-\nd,\ycoord) {};
      \node(O44) at (axis cs:TSQR-\tra-\nd,\ycoord) {};
      \node(O4) at ($(O41)!0.5!(O44)$,\ycoord) {};
      \node(N41) at (axis cs:TSQR-\trb-\nd,\ycoord) {};
      \node(N44) at (axis cs:TSQR-\trb-\nd,\ycoord) {};
      \node(N4) at ($(N41)!0.5!(N44)$,\ycoord) {};
      \node(XTL) at (xticklabel cs:0)  {};
      \node(PG1) at ($(O1)!0.5!(N1)$,\ycoord) {};
      \node(PG2) at ($(O2)!0.5!(N2)$,\ycoord) {};
      \node(PG3) at ($(O3)!0.5!(N3)$,\ycoord) {};
      \node(PG4) at ($(O4)!0.5!(N4)$,\ycoord) {};
      \node[yshift=-10,anchor=north] at (PG1 |- XTL) {\small \nalabel};
      \node[yshift=-10,anchor=north] at (PG2 |- XTL) {\small \nblabel};
      \node[yshift=-10,anchor=north] at (PG3 |- XTL) {\small \nclabel};
      \node[yshift=-10,anchor=north] at (PG4 |- XTL) {\small \ndlabel};
      \iflegend
      \else
          \node[anchor=south] at (axis cs:TSQR-\tra-\nc,151) {$\vdots$};
      \fi
    }]
    
    \addplot[color=color1,fill=color1]  table[x=alg-tr-nc, y expr=(\thisrow{TSQR_comp}/\numiters)] {\datafile};
    \addplot[color=color1,pattern=north east lines,pattern color=color1] table[x=alg-tr-nc, y expr=(\thisrow{TSQR_comm}/\numiters)] {\datafile};
    \addplot[color=color5,pattern=north east lines,pattern color=color5] table[x=alg-tr-nc, y expr=(\thisrow{other_comm}/\numiters)] {\datafile};
    \addplot[color=color2,fill=color2] table[x=alg-tr-nc, y expr=(\thisrow{AppQ_comp}/\numiters))] {\datafile};
    \addplot[color=color2,pattern=north east lines,pattern color=color2] table[x=alg-tr-nc, y expr=(\thisrow{AppQ_comm}/\numiters)] {\datafile};
    \iflegend
    	\legend{Comp QR, Comm QR, Comm R, Comp Apply-Q, Comm Apply-Q};
    \fi
  \end{axis}
}

\newcommand{\roundingplotoptions}{
  ybar stacked,
  reverse legend,
  bar width=20pt,
  \ifylabel
  	ylabel={Time (seconds)}, 
  \fi
  ymin=0,
  ymax=\na,
  symbolic x coords={LRLI, LRL, RLRI, RLR},
  xticklabels={LRLI, LRL, RLRI, RLR},  
  xtick=data,
  legend pos={north west},
  legend style={draw=none, cells={align=left}, nodes={scale=0.8}},
  minor tick num=9,
      yminorgrids,
      ymajorgrids,
}

\newcommand{\makeroundingplotbreaktime}{
  \begin{axis}[\roundingplotoptions]
    
    \addplot[color=color1,fill=color1]  table[x=alg, y expr=(\thisrow{TSQR_comp})] {\datafile};
    \addplot[color=color1,pattern=north east lines,pattern color=color1] table[x=alg-tr-nc, y expr=(\thisrow{TSQR_comm}] {\datafile};
    \addplot[color=color5,pattern=north east lines,pattern color=color5] table[x=alg-tr-nc, y expr=(\thisrow{other_comm}] {\datafile};
    \addplot[color=color2,fill=color2] table[x=alg-tr-nc, y expr=(\thisrow{AppQ_comp})] {\datafile};
    \addplot[color=color2,pattern=north east lines,pattern color=color2] table[x=alg-tr-nc, y expr=(\thisrow{AppQ_comm}] {\datafile};
    \iflegend
    	\legend{Comp QR, Comm QR, Bcast R, Comp Apply-Q, Comm Apply-Q};
    \fi
  \end{axis}
}
\newcommand{\makeroundingplotnobreaktime}{
  \begin{axis}[\roundingplotoptions]
    \addplot[color=color1,fill=color1]  table[x=alg, y expr=(\thisrow{Total})] {\datafile};
  \end{axis}
}
\newcommand{\normplotoptions}{
  ybar stacked,
  reverse legend,
  bar width=5pt,
  width = 13cm,
  	ylabel={Time fraction}, 
  ymin=0,
  ymax=1,
  symbolic x coords={\na-16,\nb-16,\nc-16,,\na-32,\nb-32,\nc-32,,\na-64,\nb-64,\nc-64,,\na-128,\nb-128,\nc-128,,\na-256,\nb-256,\nc-256},
  xticklabels={InnPro,InnPro,InnPro,InnPro,InnPro,InnPro-Sym,InnPro-Sym,InnPro-Sym,InnPro-Sym,InnPro-Sym,Ortho,Ortho,Ortho,Ortho,Ortho},  
  xtick=data,
  x tick label style={rotate=270},
  legend pos={south west},
  legend style={cells={align=left}, nodes={scale=1.6}},
}

\newcommand{\makenormplotbreaktime}{
  \begin{axis}[
      \normplotoptions
    after end axis/.code={ 
      \node(O11) at (axis cs:\na-16,1) {};
      \node(O14) at (axis cs:\na-16,1) {};
      \node(O1) at ($(O11)!0.5!(O14)$,1) {};
      \node(N11) at (axis cs:\nc-16,1) {};
      \node(N14) at (axis cs:\nc-16,1) {};
      \node(N1) at ($(N11)!0.5!(N14)$,1) {};
      \node(O21) at (axis cs:\na-32,1) {};
      \node(O24) at (axis cs:\na-32,1) {};
      \node(O2) at ($(O21)!0.5!(O24)$,1) {};
      \node(N21) at (axis cs:\nc-32,1) {};
      \node(N24) at (axis cs:\nc-32,1) {};
      \node(N2) at ($(N21)!0.5!(N24)$,1) {};
      \node(O31) at (axis cs:\na-64,1) {};
      \node(O34) at (axis cs:\na-64,1) {};
      \node(O3) at ($(O31)!0.5!(O34)$,1) {};
      \node(N31) at (axis cs:\nc-64,1) {};
      \node(N34) at (axis cs:\nc-64,1) {};
      \node(N3) at ($(N31)!0.5!(N34)$,1) {};
      \node(O41) at (axis cs:\na-128,1) {};
      \node(O44) at (axis cs:\na-128,1) {};
      \node(O4) at ($(O41)!0.5!(O44)$,1) {};
      \node(N41) at (axis cs:\nc-128,1) {};
      \node(N44) at (axis cs:\nc-128,1) {};
      \node(N4) at ($(N41)!0.5!(N44)$,1) {};
      \node(O51) at (axis cs:\na-256,1) {};
      \node(O54) at (axis cs:\na-256,1) {};
      \node(O5) at ($(O51)!0.5!(O54)$,1) {};
      \node(N51) at (axis cs:\nc-256,1) {};
      \node(N54) at (axis cs:\nc-256,1) {};
      \node(N5) at ($(N51)!0.5!(N54)$,1) {};
      \node(XTL) at (xticklabel cs:0)  {};
      \node(PG1) at ($(O1)!0.5!(N1)$,1) {};
      \node(PG2) at ($(O2)!0.5!(N2)$,1) {};
      \node(PG3) at ($(O3)!0.5!(N3)$,1) {};
      \node(PG4) at ($(O4)!0.5!(N4)$,1) {};
      \node(PG5) at ($(O5)!0.5!(N5)$,1) {};
      \node[yshift=-10,anchor=north] at (PG1 |- XTL) {\small \nalabel};
      \node[yshift=-10,anchor=north] at (PG2 |- XTL) {\small \nblabel};
      \node[yshift=-10,anchor=north] at (PG3 |- XTL) {\small \nclabel};
      \node[yshift=-10,anchor=north] at (PG4 |- XTL) {\small \ndlabel};
      \node[yshift=-10,anchor=north] at (PG5 |- XTL) {\small \nelabel};
      \iflegend
      \else
          \node[anchor=south] at (axis cs:TSQR-\tra-\nc,151) {$\vdots$};
      \fi
    }]
    
    \addplot[color=color1,fill=color1]  table[x=alg, y expr=(\thisrow{comp}/(\thisrow{comp}+\thisrow{comm}))] {\datafile};
    \addplot[color=color2,fill=color2]  table[x=alg, y expr=(\thisrow{comm}/(\thisrow{comp}+\thisrow{comm}))] {\datafile};
    \iflegend
    	\legend{Comp, Comm};
    \fi
  \end{axis}
}

\newcommand{\onelargemodenormplotoptions}{
  ybar stacked,
  reverse legend,
  bar width=5pt,
  width = 13cm,
  	ylabel={Time fraction}, 
  ymin=0,
  ymax=1,
  symbolic x coords={\na-1,,\nb-1,,\nc-1,,,\na-2,,\nb-2,,\nc-2,,,\na-4,,\nb-4,,\nc-4,,,\na-8,,\nb-8,,\nc-8,,,\na-16,,\nb-16,,\nc-16,,,\na-32,,\nb-32,,\nc-32,,,\na-64,,\nb-64,,\nc-64,,,\na-128,,\nb-128,,\nc-128,,,\na-256,,\nb-256,,\nc-256},
  xticklabels={InnPro,InnPro,InnPro,InnPro,InnPro,InnPro,InnPro,InnPro,InnPro,InnPro-Sym,InnPro-Sym,InnPro-Sym,InnPro-Sym,InnPro-Sym,InnPro-Sym,InnPro-Sym,InnPro-Sym,InnPro-Sym,Ortho,Ortho,Ortho,Ortho,Ortho,Ortho,Ortho,Ortho,Ortho},  
  xtick=data,
  x tick label style={rotate=270},
  legend pos={south west},
  legend style={cells={align=left}, nodes={scale=1.6}},
}

\newcommand{\makeonelargemodenormplotbreaktime}{
  \begin{axis}[
      \onelargemodenormplotoptions
    after end axis/.code={ 
      \node(O61) at (axis cs:\na-1,1) {};
      \node(O64) at (axis cs:\na-1,1) {};
      \node(O6) at ($(O61)!0.5!(O64)$,1) {};
      \node(N61) at (axis cs:\nc-1,1) {};
      \node(N64) at (axis cs:\nc-1,1) {};
      \node(N6) at ($(N61)!0.5!(N64)$,1) {};
      \node(O71) at (axis cs:\na-2,1) {};
      \node(O74) at (axis cs:\na-2,1) {};
      \node(O7) at ($(O71)!0.5!(O74)$,1) {};
      \node(N71) at (axis cs:\nc-2,1) {};
      \node(N74) at (axis cs:\nc-2,1) {};
      \node(N7) at ($(N71)!0.5!(N74)$,1) {};
      \node(O81) at (axis cs:\na-4,1) {};
      \node(O84) at (axis cs:\na-4,1) {};
      \node(O8) at ($(O81)!0.5!(O84)$,1) {};
      \node(N81) at (axis cs:\nc-4,1) {};
      \node(N84) at (axis cs:\nc-4,1) {};
      \node(N8) at ($(N81)!0.5!(N84)$,1) {};
      \node(O91) at (axis cs:\na-8,1) {};
      \node(O94) at (axis cs:\na-8,1) {};
      \node(O9) at ($(O91)!0.5!(O94)$,1) {};
      \node(N91) at (axis cs:\nc-8,1) {};
      \node(N94) at (axis cs:\nc-8,1) {};
      \node(N9) at ($(N91)!0.5!(N94)$,1) {};
      \node(O11) at (axis cs:\na-16,1) {};
      \node(O14) at (axis cs:\na-16,1) {};
      \node(O1) at ($(O11)!0.5!(O14)$,1) {};
      \node(N11) at (axis cs:\nc-16,1) {};
      \node(N14) at (axis cs:\nc-16,1) {};
      \node(N1) at ($(N11)!0.5!(N14)$,1) {};
      \node(O21) at (axis cs:\na-32,1) {};
      \node(O24) at (axis cs:\na-32,1) {};
      \node(O2) at ($(O21)!0.5!(O24)$,1) {};
      \node(N21) at (axis cs:\nc-32,1) {};
      \node(N24) at (axis cs:\nc-32,1) {};
      \node(N2) at ($(N21)!0.5!(N24)$,1) {};
      \node(O31) at (axis cs:\na-64,1) {};
      \node(O34) at (axis cs:\na-64,1) {};
      \node(O3) at ($(O31)!0.5!(O34)$,1) {};
      \node(N31) at (axis cs:\nc-64,1) {};
      \node(N34) at (axis cs:\nc-64,1) {};
      \node(N3) at ($(N31)!0.5!(N34)$,1) {};
      \node(O41) at (axis cs:\na-128,1) {};
      \node(O44) at (axis cs:\na-128,1) {};
      \node(O4) at ($(O41)!0.5!(O44)$,1) {};
      \node(N41) at (axis cs:\nc-128,1) {};
      \node(N44) at (axis cs:\nc-128,1) {};
      \node(N4) at ($(N41)!0.5!(N44)$,1) {};
      \node(O51) at (axis cs:\na-256,1) {};
      \node(O54) at (axis cs:\na-256,1) {};
      \node(O5) at ($(O51)!0.5!(O54)$,1) {};
      \node(N51) at (axis cs:\nc-256,1) {};
      \node(N54) at (axis cs:\nc-256,1) {};
      \node(N5) at ($(N51)!0.5!(N54)$,1) {};
      \node(XTL) at (xticklabel cs:0)  {};
      \node(PG1) at ($(O1)!0.5!(N1)$,1) {};
      \node(PG2) at ($(O2)!0.5!(N2)$,1) {};
      \node(PG3) at ($(O3)!0.5!(N3)$,1) {};
      \node(PG4) at ($(O4)!0.5!(N4)$,1) {};
      \node(PG5) at ($(O5)!0.5!(N5)$,1) {};
      \node(PG6) at ($(O6)!0.5!(N6)$,1) {};
      \node(PG7) at ($(O7)!0.5!(N7)$,1) {};
      \node(PG8) at ($(O8)!0.5!(N8)$,1) {};
      \node(PG9) at ($(O9)!0.5!(N9)$,1) {};
      \node[yshift=-10,anchor=north] at (PG1 |- XTL) {\small \nalabel};
      \node[yshift=-10,anchor=north] at (PG2 |- XTL) {\small \nblabel};
      \node[yshift=-10,anchor=north] at (PG3 |- XTL) {\small \nclabel};
      \node[yshift=-10,anchor=north] at (PG4 |- XTL) {\small \ndlabel};
      \node[yshift=-10,anchor=north] at (PG5 |- XTL) {\small \nelabel};
      \node[yshift=-10,anchor=north] at (PG6 |- XTL) {\small \nflabel};
      \node[yshift=-10,anchor=north] at (PG7 |- XTL) {\small \nglabel};
      \node[yshift=-10,anchor=north] at (PG8 |- XTL) {\small \nhlabel};
      \node[yshift=-10,anchor=north] at (PG9 |- XTL) {\small \nilabel};
      \iflegend
      \else
          \node[anchor=south] at (axis cs:TSQR-\tra-\nc,151) {$\vdots$};
      \fi
    }]
    
    \addplot[color=color1,fill=color1]  table[x=alg, y expr=(\thisrow{comp}/(\thisrow{comp}+\thisrow{comm}))] {\datafile};
    \addplot[color=color2,fill=color2]  table[x=alg, y expr=(\thisrow{comm}/(\thisrow{comp}+\thisrow{comm}))] {\datafile};
    \iflegend
    	\legend{Comp, Comm};
    \fi
  \end{axis}
}
\newcommand{\largemodesplotoptions}{
  ybar stacked,
  reverse legend,
  bar width=10pt,
  	ylabel={Time fraction}, 
  ymin=0,
  ymax=1,
  symbolic x coords={16,,32,,64,,128,,256},
  xticklabels={16,32,64,128,256},  
  xtick=data,
  legend pos={south west},
  legend style={cells={align=left}, nodes={scale=0.8}},
}

\newcommand{\makeroundinglargemodesplotbreaktime}{
  \begin{axis}[
      \largemodesplotoptions
    after end axis/.code={ 
      \iflegend
      \else
          \node[anchor=south] at (axis cs:TSQR-\tra-\nc,151) {$\vdots$};
      \fi
    }]
    
    \addplot[color=color1,fill=color1]  table[x=p, y expr=(\thisrow{TSQR_comp}/(\thisrow{Total}))] {\datafile};
    \addplot[color=color1,pattern=north east lines,pattern color=color1]  table[x=p, y expr=((\thisrow{TSQR_comm}+\thisrow{other_comm})/(\thisrow{Total}))] {\datafile};
    \addplot[color=color2,fill=color2]  table[x=p, y expr=(\thisrow{AppQ_comp}/(\thisrow{Total}))] {\datafile};
    \addplot[color=color2,pattern=north east lines,pattern color=color2]  table[x=p, y expr=(\thisrow{AppQ_comm}/(\thisrow{Total}))] {\datafile};
    \addplot[color=color3,fill=color3]  table[x=p, y expr=(\thisrow{other_comp}/(\thisrow{Total}))] {\datafile};
    \iflegend
    	\legend{TSQR Comp, TSQR Comm, AppQ Comp, AppQ Comm, Other Comp};
    \fi
  \end{axis}
}

\newcommand{\onelargemodeplotoptions}{
  ybar stacked,
  reverse legend,
  bar width=10pt,
  	ylabel={Time fraction}, 
  ymin=0,
  ymax=1,
  symbolic x coords={1,,2,,4,,8,,16,,32,,64,,128,,256},
  xticklabels={1,2,4,8,16,32,64,128,256},  
  xtick=data,
  legend pos={south west},
  legend style={cells={align=left}, nodes={scale=0.8}},
}
\newcommand{\makeroundingonelargemodeplotbreaktime}{
  \begin{axis}[
      \onelargemodeplotoptions
    after end axis/.code={ 
      \iflegend
      \else
          \node[anchor=south] at (axis cs:TSQR-\tra-\nc,151) {$\vdots$};
      \fi
    }]
    
    \addplot[color=color1,fill=color1]  table[x=p, y expr=(\thisrow{TSQR_comp}/(\thisrow{Total}))] {\datafile};
    \addplot[color=color1,pattern=north east lines,pattern color=color1]  table[x=p, y expr=((\thisrow{TSQR_comm}+\thisrow{other_comm})/\thisrow{Total})] {\datafile};
    \addplot[color=color2,fill=color2]  table[x=p, y expr=(\thisrow{AppQ_comp}/(\thisrow{Total}))] {\datafile};
    \addplot[color=color2,pattern=north east lines,pattern color=color2]  table[x=p, y expr=(\thisrow{AppQ_comm}/(\thisrow{Total}))] {\datafile};
    \addplot[color=color3,fill=color3]  table[x=p, y expr=(\thisrow{other_comp}/(\thisrow{Total}))] {\datafile};
    \iflegend
    	\legend{TSQR Comp, TSQR Comm, AppQ Comp, AppQ Comm, Other Comp};
    \fi
  \end{axis}
}

\newcommand{\makestrongscalingplot}{
\begin{axis}[
	legend pos = north east,
	reverse legend,
	ylabel={Time (seconds)}, 
	xlabel={Number of Nodes},
	xtick=data,
	xticklabels={16,32,64,128,256},
	xmode=log,
	log basis x ={2},
	ymode=log,
	log basis y={2},
	y tick label style={rotate=90},
        ylabel near ticks,
  minor tick num=9,
      grid=both,
]
	\legend{LRLI, Perfect}
	
	\addplot[color=color1, thick] table[x=p, y=runtime] {\datafile};
  \addplot[color=black, thick, dashed] table[x=p, y=perfect] {\datafile};
	

\end{axis}
}

\newcommand{\makestrongscalingplotonelargemode}{
\begin{axis}[
	legend pos = north east,
	reverse legend,
	ylabel={Time (seconds)}, 
	xlabel={Number of Nodes},
	xtick=data,
	xticklabels={1,2,4,8,16,32,64,128,256},
	xmode=log,
	log basis x ={2},
	ymode=log,
	log basis y={2},
	y tick label style={rotate=90},
        ylabel near ticks,
  minor tick num=9,
      grid=both,
]
	\legend{LRLI, Perfect}
	
	\addplot[color=color1, thick] table[x=p, y=runtime] {\datafile};
  \addplot[color=black, thick, dashed] table[x=p, y=perfect] {\datafile};
	

\end{axis}
}

\newcommand{\makestrongscalingmultiplot}{
\begin{axis}[
	legend pos = north east,
	ylabel={Time (seconds)}, 
	xlabel={Number of Nodes},
	xtick=data,
	xticklabels={16,32,64,128,256},
	xmode=log,
	log basis x ={2},
	ymode=log,
	log basis y={2},
	y tick label style={rotate=90},
        ylabel near ticks,
  minor tick num=9,
      grid=both,
]
	\legend{Perfect, Ortho, InnPro-Sym, InnPro}
	
	\addplot[color=black,thick,dashed] table[x=p, y expr=\thisrow{perfect}/\numiters] {\datafile};
	\addplot[color=color1,thick] table[x=p, y expr=\thisrow{Ortho}/\numiters] {\datafile};
	\addplot[color=color2,thick] table[x=p, y expr=\thisrow{Sym}/\numiters] {\datafile};
	\addplot[color=color3,thick] table[x=p, y expr=\thisrow{NonSym}/\numiters] {\datafile};


\end{axis}
}
\newcommand{\makeonelargemodestrongscalingmultiplot}{
\begin{axis}[
	legend pos = north east,
	ylabel={Time (seconds)}, 
	xlabel={Number of Nodes},
	xtick=data,
	xticklabels={1,2,4,8,16,32,64,128,256},
	xmode=log,
	log basis x ={2},
	ymode=log,
	log basis y={2},
	y tick label style={rotate=90},
        ylabel near ticks,
  minor tick num=9,
      grid=both,
]
	\legend{Perfect, Ortho, InnPro-Sym, InnPro}
	
	\addplot[color=black,thick,dashed] table[x=p, y expr=\thisrow{perfect}/\numiters] {\datafile};
	\addplot[color=color1,thick] table[x=p, y expr=\thisrow{Ortho}/\numiters] {\datafile};
	\addplot[color=color2,thick] table[x=p, y expr=\thisrow{Sym}/\numiters] {\datafile};
	\addplot[color=color3,thick] table[x=p, y expr=\thisrow{NonSym}/\numiters] {\datafile};


\end{axis}
}

In this section we present performance results for TT computations using synthetic tensors with mode and dimension parameters inspired by physics and chemistry applications, as described in \Cref{sec:perf_syn}.
We first present microbenchmarks in \Cref{sec:perf_micro} to justify key design decisions, and then demonstrate performance efficiency and parallel scaling in \Cref{sec:perf_scaling}. 

All numerical experiments are run on the Max Planck Society supercomputer COBRA.
All computation nodes contain two Intel Xeon Gold 6148 processors (Skylake, 20 cores each at 2.4 GHz) and 192 GB of memory, and the nodes are connected through a 100 Gb/s OmniPath interconnect.
We link to MKL 2020.1 for single-threaded BLAS and LAPACK subroutines.

\subsection{\change{Motivating Applications}}

\change{We describe in this section the motivating applications guiding the choice of tensor dimensions and ranks of the synthetic models we use in the experiments.}

\subsubsection{\change{High-Order Correlation Functions}}\label{sec:HOCF}
In the study of stochastic processes, Gaussian random fields are widely used.
If $f$ is a Gaussian random field defined on a bounded domain $\Omega \subset \R^d$ ($d=1,2,3$), an $N$-point correlation function for $f$ is defined on $\Omega^N$.
\change{The discretization of the domain determines the $N$-way tensor dimensions.}
These $N$-point correlation functions can often be efficiently approximated in TT format \cite{BonNK16,KreKNT15}.
\change{For typical discretizations, the number of discretization points in the domain $\Omega$ can be extremely large leading to tensors with very large dimensions.
In order to compute some desired information about the random solution of a stochastic PDE such as its expected value, TT computations including addition and scaling are required.}
Thus, compressing the resulting TT tensors is required to maintain the tractability of computations.
\change{In \cite{BonNK16} the authors present a study of single-phase fluid flows in heterogeneous porous media. Due to memory and time constraints, current implementations of TT arithmetic allows only to perform the aforementioned computations on a
moderate size discretizations (10,000) for $d=1$ or $d=2$. However, in industrial applications where $\Omega \subset \R^3$, the mode dimension can be of order $10^8$.}

\subsubsection{\change{Molecular Simulations}}\label{sec:MS}
Another important class of applications is molecular simulations.
For example, when a spin system can be considered as a weakly branched linear chain, it is typical to represent it as a TT tensor \cite{SavDWK14}.
Each branch is then considered as a spatial coordinate (mode).
The number of branches\change{, corresponding to the number of tensor modes,} can be arbitrarily large; for example, a simple backbone protein may have hundreds of branches.
The TT representation is then inherited from the weak correlation between the branches.
However, in the same branch, the correlation cannot be ignored, and thus the exponential growth in the number of states\change{, which corresponds to the dimension of the tensor mode for that branch,} cannot be avoided.

\subsubsection{\change{Parameter-Dependent PDEs}}\label{sec:PDPDES}
In this application, one or a few modes may be much larger than the rest.
This is typically the case in physical applications such as parameter-dependent PDEs, stochastic PDEs, uncertainty quantification, and optimal control systems \cite{BenDOS16,BenDOS20,BenGW15,DolS17,HesRS15,KreT10,QuaMN15}.
In such applications, the spatial discretization leads to a high number of degrees of freedom.
This typically results from large domains, refinement procedures, and a large number of parameter samples.
Most of other modes correspond to control or uncertainty parameters and can have relatively smaller dimension.

\change{For example, in \cite{BenDOS20} where the authors study an optimal control problem constrained by random Navier--Stokes equations, certain vectors are represented by 10-mode tensors. The number of degrees of freedom in each mode is as follows:
the velocity field has up to 168,240, the time mode has up to 4096, and the eight modes related to the random variables each has 8.}
\change{Again, this discretization is limited by memory and time constraints and finer granularity that increases the accuracy of the approximation would lead to dimensions on the order of millions.}
\subsection{Synthetic TT Models}
\label{sec:perf_syn}

As we are interested in large scale systems, we consider two contexts of applications in which a large number of modes exists.
The first context is with each mode of relatively the same (large) dimension\change{, such as the applications described in \Cref{sec:HOCF,sec:MS}}, and the second context is a single or few modes with large dimension as well as many modes of relatively smaller dimension\change{, as arises in parameter-dependent PDEs (\Cref{sec:PDPDES})}.
\Cref{tab:models} presents the details of the three models of synthetic tensors we use in the experiments, in order of their memory size.
The first and third models correspond to the first context (all modes of the same dimension) and the second model corresponds to the second context (two large modes and many more smaller modes).
The first model is chosen to be small enough to be processed by a single core, while the second and third are larger and benefit more from distributed-memory parallelization (the third does not fit in the memory of a single node).
The paragraphs below describe the applications that inspire these choices of modes and dimensions.

In all experiments, we generate a random TT tensor $\T{X}$ with a given number of modes $N$, modes sizes $I_n$ for $n = 1, \ldots, N$, and TT ranks $R^{\T{X}}_n$ for $n = 1, \ldots, N - 1$.
Then, we form the TT tensor $ \T{Y} = 2 \T{X} - \T{X}$ \change{whose representation has TT} ranks $R^{\T{Y}}_n = 2 R^{\T{X}}_n$ for $n = 1, \ldots, N - 1$.
The algorithms are then applied on the TT tensor $\T{Y}$.
Note that the minimal TT ranks of $\T Y$ are less or equal than the TT ranks of $\T X$.

\begin{table}
\centering
\begin{tabular}{c|c|c|c|c}
 Model & \# Modes & Dimensions & Ranks & Memory \\
\hline
 1 & 50 & $2K \times \cdots \times 2K$ & $50$ & 2 GB \\ 
 2 & 16 & $100M \times 50K \times \cdots \times 50K \times 1M$ & $30$ & 28 GB \\ 
 3 & 30 & $2M \times \cdots \times 2M$ & $30$ & 385 GB \\ 
\hline
\end{tabular}
\caption{Synthetic TT models used for performance experiments.  In each case the \change{TT} ranks are all the same and are cut in half by the TT rounding procedure.}
\label{tab:models}
\end{table}

\subsection{Microbenchmarks}
\label{sec:perf_micro}

We next present experimental results for microbenchmarks to justify our choices for subroutine algorithms and optimizations.
The results presented in \Cref{sec:perf_scaling} use the best-performing variants and optimizations demonstrated in this section.

\subsubsection{TSQR}
\label{sec:perf_tsqr}

As discussed in \Cref{sec:TSQR}, the TSQR algorithm depends on a hierarchical tree.
Two tree choices are commonly used in practice, the binomial tree and the butterfly tree.
In both cases the TSQR computes the QR decomposition sharing the same complexity and communication costs along the critical path, whereas the butterfly requires less communication cost along the critical path of the application of the implicit \change{orthonormal} factor.
\change{This advantage of the butterfly variant in the application phase is particularly important in the context of TT orthonormalization and rounding because a large percentage of time is spent in the application phase.}

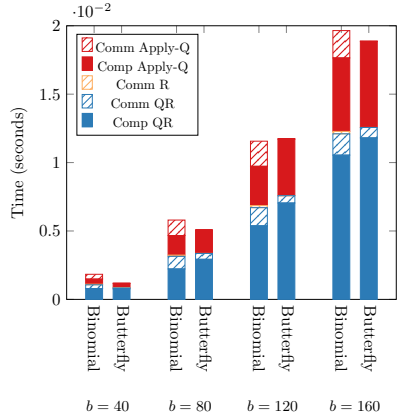
\begin{figure}
  \centering
  \begin{adjustbox}{width=0.4\linewidth}
  \begin{tikzpicture}[scale=1]
    \renewcommand{\datafile}{data/tsqr.dat}
    \renewcommand{\na}{40}
    \renewcommand{\nb}{80}
    \renewcommand{\nc}{120}
    \renewcommand{\nd}{160}
    \renewcommand{\nalabel}{$b=40$}
    \renewcommand{\nblabel}{$b=80$}
    \renewcommand{\nclabel}{$b=120$}
    \renewcommand{\ndlabel}{$b=160$}
    \renewcommand{\tra}{BN}
    \renewcommand{\trb}{BF}
    \renewcommand{\numiters}{50}
    \legendtrue
    \ylabeltrue
    \maketsqrplot
  \end{tikzpicture}
  \end{adjustbox}
  \caption{Time breakdown for TSQR variants for $1{,}024{,}000 \times b$ matrix over 1024 processors, including both factorization and application of the \change{orthonormal} factor to a dense $b\times b$ matrix.}
  \label{fig:tsqr}
\end{figure}

Here we compare the performance of the TSQR algorithms using the binomial and butterfly trees for both factorization and single application of the \change{orthonormal} factor.
Since the difference in their costs is solely related to the number of columns, we fix the number of rows in the comparison and vary the number of columns.
\Cref{fig:tsqr} reports the breakdown of time of the variants using 256 nodes with 4 MPI processes per node (2 cores per socket). 
The local matrix size on each processor is $1{,}000 \times b$ where $b$ varies in $\{40, 80, 120, 160\}$.
We observe that the butterfly tree has better performance in terms of communication time in the application phase.
Note that the factorization runtime (computation and communication) is relatively the same for both variants.
We also time the cost of communicating the triangular factor $\M R$, which is required of the binomial variant in the context of TT-rounding, but that cost is negligible in these experiments.

Based on these results (and corroborating experiments with various other parameters), we use the butterfly variant of TSQR for TT computations that require TSQR in all subsequent numerical experiments.

\subsubsection{TT Rounding}
\label{sec:micro_rounding}

In this section, we consider 4 variants of TT rounding (\Cref{alg:par_tt_rounding}), based on the \change{orthonormal}ization/truncation ordering and the use of the implicit \change{orthonormal} factor optimization.
As discussed in \Cref{sec:seq-tt-rounding}, the rounding procedure can perform right- or left- \change{orthonormal}ization followed by a truncation phase in the opposite direction.
We refer to the ordering based on right-\change{orthonormal}ization and left-truncation as RLR and the ordering based on left-\change{orthonormal}ization and right-truncation as LRL.
The implicit optimization avoids the explicit formation of \change{orthonormal} factors during the \change{orthonormal}ization phase; instead of using \Cref{alg:par_tt_orthonormalization} as a black-box subroutine, \Cref{alg:par_tt_rounding} leaves \change{orthonormal} factors in implicit TSQR form as much as possible, saving a constant factor of computation (and a small amount of communication).

Although the asymptotic complexity of the variants of the rounding procedure are equal, their performance is not the same.
This disparity between RLR and LRL orderings is because of the performance difference between the QR and the LQ implementations of the LAPACK subroutines provided by the MKL implementations.
Despite the same computation complexity, the QR subroutines has much better performance than the LQ subroutines.

In the LRL ordering, a sequence of calls to the QR subroutine are performed on the vertically unfolded TT cores $\TT{X}{n}$ with the increased ranks $R_{n-1}, R_n$.
Along the truncation sweep, the LQ subroutine is called in a sequence to factor the horizontally unfolded TT cores $\TT{X}{n}$ with one reduced rank $R_{n-1}, L_n$.
As presented in \Cref{sec:par_orth,sec:par_round}, the RLR ordering employs the QR and LQ subroutines in the opposite order. 
Because the truncation phase involves less computation within local QR/LQ subroutine calls than the \change{orthonormal}ization phase, the LRL ordering has the advantage that it spends less time in LQ subroutine calls than the RLR ordering.

The effect of the implicit optimization is a reduction in computation (approximately 12.5\% in these experiments) and communication, but this advantage is offset in part by the performance of local subroutines.
The implicit application of the \change{orthonormal} factor involves auxiliary LAPACK routines for applying sets of Householder vectors in various formats.
The explicit multiplication of an \change{orthonormal} factor to a small square matrix involves a broadcast and a local subroutine call to matrix multiplication, which has much higher performance than the auxiliary routines involving Householder vectors.
We use an ``I'' to indicate the use of the implicit optimization, so that the 4 variants are LRLI, LRL, RLRI, and RLR.

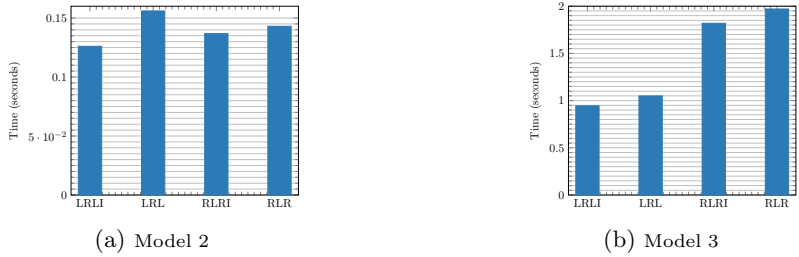
\begin{figure}
\begin{subfigure}[b]{0.5\linewidth}
  \centering
  \begin{adjustbox}{width=0.6\linewidth}
  \begin{tikzpicture}
    \renewcommand{\datafile}{data/tt_one_large_mode_1280_nobreaktime.dat}
    \legendfalse
    \ylabeltrue
    \renewcommand{\na}{0.16}
    \makeroundingplotnobreaktime
  \end{tikzpicture}
  \end{adjustbox}
  \captionsub{Model 2}
  \label{fig:tt_one_large_mode_1280_nobreaktime}
\end{subfigure}
\hfill
\begin{subfigure}[b]{0.46\linewidth}
  \centering
  \begin{adjustbox}{width=0.6\linewidth}
  \begin{tikzpicture}
    \renewcommand{\datafile}{data/tt_large_modes_1280_nobreaktime.dat}
    \legendfalse
    \ylabeltrue
    \renewcommand{\na}{2}
    \makeroundingplotnobreaktime
  \end{tikzpicture}
  \end{adjustbox}
  \captionsub{Model 3}
  \label{fig:tt_large_modes_1280_nobreaktime}
\end{subfigure}
\caption{Performance comparison of TT-Rounding variants for large TT models on 32 nodes (1,280 cores).  LRL refers to left-\change{orthonormal}ization followed by right-truncation (vice versa for RLR) and I indicates the use of the implicit optimization.}
\label{fig:rounding-micro}
\end{figure}

\Cref{fig:rounding-micro} presents the performance results for TT Models 2 and 3 running on 256 nodes.
We see that for both models, the LRL ordering with the implicit optimization (LRLI) is the fastest.
In the case of Model 2, the implicit optimization makes more of a difference than the ordering.
This is because a considerable amount of time is spent in the first mode, where the QR is used (once) in either ordering.
In the case of Model 3, the ordering makes a much larger difference in running time, as the internal modes dominate the running time and the QR/LQ difference has a large effect.
The implicit optimization still improves performance, but it has less of an effect than the ordering.
Based on these results, we use the LRLI variant of TT-rounding in all the experiments presented in \Cref{sec:perf_scaling}.

\subsection{Parallel Scaling}
\label{sec:perf_scaling}

\subsubsection{Norms}
\label{sec:normscaling}

In this section we compare the performance and parallel scaling of three different algorithms for computing the norm of a TT tensor as discussed in \Cref{sec:norms}.
We focus on this computation because the multiple approaches represent the performance of algorithms for computing inner products and \change{orthonormal}ization, which are essential on their own in other contexts.
We use ``Ortho'' to denote the approach of first right- or left-\change{orthonormal}izing the TT tensor and then (cheaply) computing the norm of the first or last core, respectively.
Thus, Ortho performance represents that of \Cref{alg:par_tt_orthonormalization}.
The name ``InnPro'' refers to the approach of computing the inner product of the TT tensor with itself, and ``InnPro-Sym'' includes the optimization that exploits the symmetry in the inner product to save up to half the computation.
InnPro captures the performance of the algorithm described in \Cref{sec:inner} for general TT inner products as well.

\begin{figure}
  \centering
    \begin{subfigure}[b]{0.49\linewidth}
    \centering
    \begin{adjustbox}{width=0.85\linewidth}
      \begin{tikzpicture}[scale=.9]
        \renewcommand{\datafile}{data/tt_one_large_mode_norm_breaktime.dat}
        \renewcommand{\na}{Had}
        \renewcommand{\nb}{Sym}
        \renewcommand{\nc}{Ort}
        \newcommand{\nflabel}{1}
        \newcommand{\nglabel}{2}
        \newcommand{\nhlabel}{4}
        \newcommand{\nilabel}{8}
        \renewcommand{\nalabel}{16}
        \renewcommand{\nblabel}{32}
        \renewcommand{\nclabel}{64}
        \renewcommand{\ndlabel}{128}
        \renewcommand{\nelabel}{256}
        \legendtrue
        \ylabeltrue
        \makeonelargemodenormplotbreaktime
      \end{tikzpicture}
    \end{adjustbox}
    \captionsub{Time breakdown for Model 2}
    \label{fig:one_large_mode_norm_breaktime}
  \end{subfigure}
    \hfill
    \begin{subfigure}[b]{0.49\linewidth}
      \vfill
      \centering
      \begin{adjustbox}{width=0.85\linewidth}
        \begin{tikzpicture}[scale=.9]
          \renewcommand{\datafile}{data/tt_large_modes_norm_breaktime.dat}
          \renewcommand{\na}{Had}
          \renewcommand{\nb}{Sym}
          \renewcommand{\nc}{Ort}
          \renewcommand{\nalabel}{16}
          \renewcommand{\nblabel}{32}
          \renewcommand{\nclabel}{64}
          \renewcommand{\ndlabel}{128}
          \renewcommand{\nelabel}{256}
          \legendtrue
          \ylabeltrue
          \makenormplotbreaktime
        \end{tikzpicture}
      \end{adjustbox}
      \captionsub{Time breakdown for Model 3}
      \label{fig:large_modes_norm_breaktime}
    \end{subfigure}
    \\
      \begin{subfigure}[b]{0.49\linewidth}
    \centering
    \begin{adjustbox}{width=0.85\linewidth}
      \begin{tikzpicture}[scale=.9]
        \renewcommand{\numiters}{5}
        \renewcommand{\datafile}{data/tt_one_large_mode_norm_scaling.dat}
        \makeonelargemodestrongscalingmultiplot
      \end{tikzpicture}
    \end{adjustbox}
    \captionsub{Parallel scaling for Model 2}
    \label{fig:one_large_mode_norm_scaling}
  \end{subfigure}
    \begin{subfigure}[b]{0.49\linewidth}
      \vfill
      \centering
      \begin{adjustbox}{width=0.85\linewidth}
        \begin{tikzpicture}[scale=.9]
          \renewcommand{\numiters}{5}
          \renewcommand{\datafile}{data/tt_large_modes_norm_scaling.dat}
          \makestrongscalingmultiplot
        \end{tikzpicture}
      \end{adjustbox}
      \captionsub{Parallel scaling for Model 3}
      \label{fig:large_modes_norm_scaling}
    \end{subfigure}
  \\
  \caption{Time breakdown and and parallel scaling of variants for TT norm computation.  ``Ortho'' refers to \change{orthonormal}ization (following by computing the norm of a single core), ``InnPro" refers to using the inner product algorithm, and ``InnPro-Sym'' refers to using the inner product algorithm with symmetric optimization.}
  \label{fig:norm}
\end{figure}

We report parallel scaling and a breakdown of computation and communication for all three algorithms and TT Models 2 and 3 in \Cref{fig:norm}.
Model 2 can be processed on a single node, but Model 3 requires 16 nodes to achieve sufficient memory; we scale both models up to 256 nodes (10,240 cores).
Based on the theoretical analysis (see \Cref{tab:costs}), when all tensor dimensions are equivalent such as Model 3, Ortho has a leading-order flop constant of 5,  InnPro has a constant of 4, and InnPro-Sym has a constant of 2.
Ortho also requires more complicated TSQR reductions compared to the All-Reduces performed in InnPro and InnPro-Sym, involving an extra $\log P$ factor in data communicated in theory and slightly less efficient implementations in practice.
In addition, the efficiencies of the local computations differ across approaches: Ortho is bottlenecked by local QR, InnPro by local matrix multiplication (GEMM), and InnPro-Sym by local triangular matrix multiplication (TRMM).

Overall, we see that InnPro is typically the best performing approach.
The main factor in its superiority is that its computation is cast as GEMM calls, which are more efficient than TRMM and QR subroutines.
Although InnPro-Sym performs half the flops of InnPro, the relative inefficiency of those flops translates to a less than $2\times$ speedup over InnPro for Model 3 and a slight slowdown for Model 2.
We also note that for high node counts, the cost of the LDLT factorization performed within InnPro-Sym becomes nonneglible and begins to hinder parallel scaling.

Based on the breakdown of computation and communication, we see that all three approaches are able to scale reasonably well because they remain computation bound up to 256 nodes.
For Model 2, we see that communication costs are relatively higher, as that tensor is much smaller.
Note that Ortho scales better than InnPro-Sym and InnPro, even superlinearly for Model 3, which is due in large part to the higher flop count and relative inefficiency of the local QRs, allowing it to remain more computation bound than the alternatives.
Overall, these results confirm that the parallel distribution of TT cores allows for high performance and scalability of the basic TT operations as described in \Cref{sec:basicops}.

\subsubsection{TT Rounding}

\paragraph{Single-Node Performance}

We compare in this section our implementation of TT rounding against the MATLAB TT-Toolbox \cite{TTToolbox} rounding process.
\Cref{tab:singlenode} presents a performance comparison on a single node of COBRA, which has 40 cores available.
We run the experiment on TT Model 1, which is small enough to be processed by a single core.
Because it is written in MATLAB, the TT-Toolbox accesses the available parallelism only through underlying calls to a multithreaded implementation of BLAS and LAPACK.
However, the bulk of the computation occurs in MATLAB functions that make direct calls to efficient BLAS and LAPACK subroutines, so it can achieve relatively high sequential performance.

We observe from \Cref{tab:singlenode} that the single-core performance of the two implementations is similar, with a 70\% speedup from our implementation.
The single-core implementations are employing the same algorithm, and we attribute the speedup to our lower-level interface to LAPACK subroutines and the ability to maintain implicit \change{orthonormal} factors to reduce computation.
The parallel strong scaling differs more drastically, as expected.
The MATLAB implementation, which is not designed for parallelization, achieves less than a $2\times$ speedup when using 20 or 40 cores.
Our parallelization, which is designed for distributed-memory systems, also scales very well on this shared-memory machine, achieving over $20\times$ speedup on 20 cores and $34\times$ speedup on 40 cores.

\begin{table}
\centering
  \begin{adjustbox}{width=\linewidth}
    \begin{tabular}{c||c||c|c||c|c}
                       & 1 core               & 20 cores              & Par.~Speedup          & 40 cores              & Par.~Speedup \\
     \hline
    TT-Toolbox         & 15.68                & 8.34                  & $\mathbf{1.9\times}$  & 8.752                 & $\mathbf{1.8\times}$ \\
    Our Implementation & 9.2                  & 0.44                  & $\mathbf{20.9\times}$ & 0.27                  & $\mathbf{33.9\times}$ \\
    \hline
      Speedup          & $\mathbf{1.7\times}$ & $\mathbf{18.95\times}$ &                       & $\mathbf{32.2\times}$ & \\
    \end{tabular}
  \end{adjustbox}
\caption{Single-node performance results on TT Model 1 and comparison with the MATLAB TT-Toolbox.}
\label{tab:singlenode}
\end{table}

\paragraph{Distributed-Memory Strong Scaling}

\begin{figure}
    \begin{subfigure}[b]{0.49\linewidth}
    \centering
    \begin{adjustbox}{width=0.65\linewidth}
      \begin{tikzpicture}[scale=.9]
        \renewcommand{\datafile}{data/one_large_mode_breaktime_Timing_pcttBFCompressImplicit_R.dat}
        \legendtrue
        \ylabeltrue
        \makeroundingonelargemodeplotbreaktime
      \end{tikzpicture}
    \end{adjustbox}
    \captionsub{Time breakdown for Model 2}
    \label{fig:one_large_mode_rounding_breaktime}
  \end{subfigure}
    \hfill
    \begin{subfigure}[b]{0.49\linewidth}
      \vfill
      \centering
      \begin{adjustbox}{width=0.65\linewidth}
        \begin{tikzpicture}[scale=.9]
          \renewcommand{\datafile}{data/large_modes_breaktime_Timing_pcttBFCompressImplicit_R.dat}
          \legendtrue
          \ylabeltrue
          \makeroundinglargemodesplotbreaktime
        \end{tikzpicture}
      \end{adjustbox}
      \captionsub{Time breakdown for Model 3}
      \label{fig:large_modes_rounding_breaktime}
    \end{subfigure}
    \\
\begin{subfigure}[b]{0.49\linewidth}
  \centering
  \begin{adjustbox}{width=0.65\linewidth}
  \begin{tikzpicture}
    \renewcommand{\datafile}{data/tt_one_large_mode_nobreaktime_scaling.dat}
    \makestrongscalingplotonelargemode
  \end{tikzpicture}
  \end{adjustbox}
  \captionsub{Parallel scaling for Model 2}
  \label{fig:tt_one_large_mode_strong_scaling}
\end{subfigure}
\begin{subfigure}[b]{0.49\linewidth}
  \centering
  \begin{adjustbox}{width=0.65\linewidth}
  \begin{tikzpicture}
    \renewcommand{\datafile}{data/tt_large_modes_nobreaktime_scaling.dat}
    \makestrongscalingplot
  \end{tikzpicture}
  \end{adjustbox}
  \captionsub{Parallel scaling for Model 3}
   \label{fig:tt_large_modes_strong_scaling}
  \end{subfigure}
  \caption{Time breakdown and and parallel scaling of LRLI variant of TT rounding.}
  \label{fig:roundingscaling}
\end{figure}
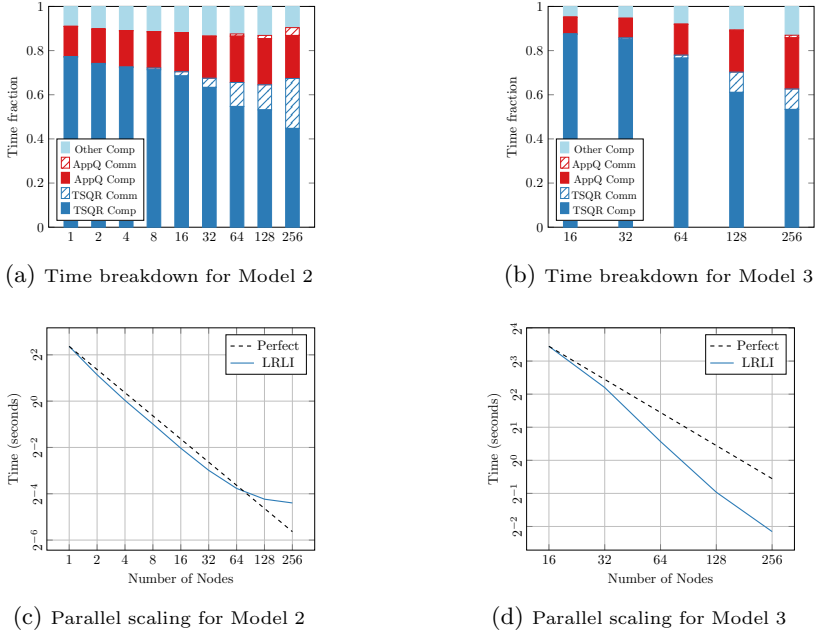

We now present the parallel performance of TT rounding scaling up to hundreds of nodes (over 10,000 cores).
As in the case of \Cref{sec:normscaling}, we consider Models 2 and 3.
\Cref{fig:roundingscaling} presents the relative time breakdown and raw timing numbers for each model.
We use the `LRLI' variant of TT rounding in these experiments per the results of \Cref{sec:micro_rounding}.
As in other rounding experiments, the ranks are cut in half for each model.

In the time breakdown plots of \Cref{fig:one_large_mode_rounding_breaktime,fig:large_modes_rounding_breaktime}, we distinguish among TSQR factorization (TSQR), application of \change{orthonormal} factors (AppQ), and the rest of the computation that includes SVDs and triangular multiplication (Other).
We also separate the computation and communication of each category.
In the context of \Cref{alg:par_tt_rounding}, TSQR corresponds to \cref{line:R2LQRi,line:L2RQRi}, AppQ corresponds to \cref{line:L2RappQi,line:L2RappVi}, and Other corresponds to \cref{line:R2LappRi,line:L2RSVDi}.

In \Cref{fig:tt_one_large_mode_strong_scaling,fig:tt_large_modes_strong_scaling}, we observe the strong scaling raw times in log scale compared to perfect scaling (based on time at the fewest number of nodes).
We see nearly perfect scaling for Model 2 until 128 nodes; time continues to decrease but is not cut in half when scaling to 256 nodes.
The parallel speedup numbers for Model 2 are $97\times$ for 128 nodes and $108\times$ for 256 nodes, compared to performance on 1 node.
In the case of Model 3, we see super-linear scaling, even at 256 nodes.
We attribute this scaling in part to the baseline comparison of 16 nodes, which already involves parallelization/communication, and in part to local data fitting into higher levels of cache as the number of processors increases, which helps memory-bound local computations.
We observe a $48\times$ speedup for Model 3, scaling from 16 to 256 nodes.

The time breakdown plots also help to explain the scaling performance.
We see that for Model 2, over 70\% of the time is spent in local computation, while for Model 3, over 90\% of the time is computation.
Of this computation, the majority is spent in TSQR, which itself is dominated by the initial local leaf QR computations.
If the rank is reduced by a smaller factor, then relatively more flops will occur in AppQ.
We note that AppQ involves minimal communication because of the use of the Butterfly TSQR variant.
The Other category is dominated by the triangular matrix multiplication, which achieves higher performance than the LAPACK subroutines involving \change{orthonormal} factors.

\section{Conclusions}
\label{sec:conclusion}

This work presents the parallel implementation of the basic computational algorithms for tensors represented in low-rank TT format.
Because most TT computations involve dependence through the train, we specify a data distribution that distributes each core across all processors and show that the computations and communication costs of our proposed algorithms enable efficiency and scalability for each core computation.
The \change{orthonormal}ization and rounding procedures for TT tensors depend heavily on the TSQR algorithm, which is designed to scale well on architectures with a large number of processors for matrices with highly skewed aspect ratios. 
Our numerical experiments show that our algorithms are indeed efficient and scalable, outperforming productivity-oriented implementations on a single core and single node and scaling well to hundreds of nodes (thousands of cores).
Thus, our approach is useful to applications and users who are restricted to a single workstation as well to those requiring the memory and performance of a supercomputer.

We note that the raw performance of our implementation depends heavily on the local BLAS/LAPACK implementation and the efficiency of the QR decomposition and related subroutines.
For example, we observe significant performance differences between MKL's implementations of QR and LQ subroutines, which caused the LRL ordering of TT-rounding to outperform RLR.
We also observe performance differences among other subroutines, such as triangular matrix multiplication and general matrix multiplication, again confirming that simple flop counting (even tracking constants closely) does not always accurately predict running times.

There do exist limitations of the parallelization approach proposed in this paper.
In particular, modes with small dimensions benefit less from parallelization and can become bottlenecks if there are too many of them.
For example, we see the limits of scalability with TT Model 2, which has large first and last modes but smaller internal modes.
In fact, the distribution scheme assumes that $P\leq I_n$ for $n= 1, \ldots, N$, and involves idle processors when the assumption is broken.
We also note that TSQR may not be the optimal algorithm to factor the unfolding, which can happen if two successive ranks differ greatly and $P$ is large with respect to the original tensor dimensions. 

Alternative possibilities to avoid these limitations include cheaper but less accurate methods for the SVD, including via the associated Gram matrices or by using randomization.
We plan to pursue such strategies in the future, in addition to considering the case of computing a TT approximation from a tensor in explicit full format.
Given these efficient computational building blocks, the next step is to build scalable Krylov and alternating-scheme based solvers that exploit the TT format.

\appendix

\section{TT Rounding Identity}
\label{app:quadprod}
We provide the full derivation of \cref{eq:quadprodTT}, which we repeat here.
The unfolding of $\T{X}$ that maps the first $n$ tensor dimensions to rows can be expressed as a product of four matrices:
\begin{equation*}
\Mz{X}{1:n} = (\M{I}_{I_n} \otimes \Mz{Q}{1:n-1}) \cdot \VOp(\TT{X}{n}) \cdot \HOp(\TT{X}{n+1}) \cdot (\M{I}_{I_{n+1}} \otimes \Mz{Z}{1}),
\end{equation*}
where $\T Q$ is $I_1 \times \cdots \times I_{n-1} \times R_{n-1}$ with 
$$\T{Q}(i_1,\dots,i_{n-1},r_{n-1}) = \TTslice{X}{1} \cdot \TTslice{X}{2} \cdots \TT{X}{n-1}(:,i_{n-1},r_{n-1}),$$
and $\T Z$ is $R_{n+1} \times I_{n+2} \times \dots \times I_N$ with
$$\T{Z}(r_{n+1},i_{n+2},\dots,i_N) = \TT{X}{n+2}(r_{n+1},i_{n+2},:) \cdot \TTslice{X}{n+3} \cdots \TTslice{X}{N}.$$
\change{The TT rounding process truncates the rank of this unfolding for each $1\leq n<N$, reducing the dimension $R_n$ to a smaller value subject to the approximation error threshold.}

Let $\T{U}$ be $I_1\times \cdots\times I_n \times R_n$ such that $\Mz{U}{1:n} = (\M{I}_{I_n} \otimes \Mz{Q}{1:n-1}) \VOp(\TT{X}{n})$, then
\begin{align*}
  \T{U}(i_1,\dots,i_n,r_n) &= \sum_{i_n'} \sum_{r_{n-1}} \delta_{(i_n',i_n)} \T{Q}(i_1,\dots,i_{n-1},r_{n-1}) \TT{X}{n}(r_{n-1},i_n',r_n) \\
 &= \sum_{r_{n-1}} \T{Q}(i_1,\dots,i_{n-1},r_{n-1}) \TT{X}{n}(r_{n-1},i_n,r_n) \\
 &= \T{Q}(i_1,\dots,i_{n-1},:) \cdot \TT{X}{n}(:,i_n,r_n) \\
 &= \TTslice{X}{1} \cdots \TTslice{X}{n-1} \cdot \TT{X}{n}(:,i_n,r_n).
\end{align*}
Let $\T{V}$ be $R_n\times I_{n+1}\times \cdots \times I_N$ such that $\Mz{V}{1} = \HOp(\TT{X}{n+1}) (\M{I}_{I_{n+1}} \otimes \Mz{Z}{1})$, then
\begin{align*}
  \T{V}(r_n,i_{n+1},\dots,i_N) &= \sum_{i_{n+1}'} \sum_{r_{n+1}} \TT{X}{n+1}(r_n,i_{n+1}',r_{n+1}) \delta_{(i_{n+1}',i_{n+1})} \T{Z}(r_{n+1},i_{n+2},\dots,i_N)  \\
 &= \sum_{r_{n-1}} \TT{X}{n+1}(r_n,i_{n+1},r_{n+1}) \T{Z}(r_{n+1},i_{n+2},\dots,i_N)  \\
 &= \TT{X}{n+1}(r_n,i_{n+1},:) \cdot \T{Z}(:,i_{n+2},\dots,i_N)  \\
 &= \TT{X}{n+1}(r_n,i_{n+1},:) \cdot \TTslice{X}{n+2} \cdots \TTslice{X}{N}.
\end{align*}

Then we confirm that $\T{Y}=\T{X}$ for $\Mz{Y}{1:n}=\Mz{U}{1:n} \cdot \Mz{V}{1}$:
\begin{align*}
\T{Y}(i_1,\dots,i_N) &=\sum_{r_n} \T{U}(i_1,\dots,i_n,r_n) \T{V}(r_n,i_{n+1},\dots,i_N) \\
 &= \sum_{r_n} \TTslice{X}{1} \cdots \TTslice{X}{n-1} \cdot \TT{X}{n}(:,i_n,r_n) \cdot\\
 & \qquad \qquad \TT{X}{n+1}(r_n,i_{n+1},:) \cdot \TTslice{X}{n+2} \cdots \TTslice{X}{N} \\
 &= \TTslice{X}{1} \cdots \TTslice{X}{n-1} \cdot\\
  & \qquad \qquad \left(\sum_{r_n} \TT{X}{n}(:,i_n,r_n) \cdot \TT{X}{n+1}(r_n,i_{n+1},:)\right) \cdot \\
 & \qquad \qquad \qquad \qquad \TTslice{X}{n+2} \cdots \TTslice{X}{N} \\
 &= \TTslice{X}{1} \cdots \TT{X}{n}(:,i_n,:) \cdot \TT{X}{n+1}(:,i_{n+1},:) \cdots \TTslice{X}{N}.
\end{align*}

\section{TSQR Subroutines for Non-Powers-of-Two}
\label{app:tsqr}
We provide here the full details of the butterfly TSQR algorithm and the algorithm for applying the resulting implicit \change{orthonormal} factor to a matrix.
These two algorithms generalize \Cref{alg:BFTSQR,alg:BFappQ} presented in \Cref{sec:TSQR} which can run only on powers-of-two processors.
To handle a non-power-of-two number of processors, we consider the first $2^{\lfloor \log P \rfloor}$ processors to be ``regular'' processors and the last $P - 2^{\lfloor \log P \rfloor}$ processors to be ``remainder'' processors.
Each remainder process has a partner in the set of regular processors, and we perform cleanup steps between remainder processors and their partners before and after the regular butterfly loop of the TSQR algorithm.
For the application algorithm, the clean up occurs after the butterfly on the regular processors (which requires no communication) and involves a single message between remainder processors and their partners.
We note that the notation and indexing matches that of \Cref{alg:BFTSQR,alg:BFappQ}, so that the algorithms coincide when $P$ is a power of two.

\begin{algorithm}
\caption{Parallel Butterfly TSQR}
\label{alg:BFTSQRapp}
\begin{algorithmic}[1]
\Require $\M{A}$ is an $m\times b$ matrix 1D-distributed so that proc $p$ owns row block $\M{A}^{(p)}$
  \Ensure $\M{A}=\M{Q} \M{R}$ with $\M{R}$ owned by all procs and $\M{Q}$ represented by $\{\M{Y}_{\ell}^{(p)}\}$ with redundancy $\M{Y}_{\ell}^{(p)}=\M{Y}_{\ell}^{(q)}$ for $p \equiv q \mod 2^\ell$ where $p, q < 2^{\lfloor \log P \rfloor}$ and $l < \lceil \log P \lceil$
\Function{$[\{\M{Y}_{\ell}^{(p)}\}, \M{R}] = $ Par-TSQR}{$\M{A}^{(p)}$}
	\State $p = \textsc{MyProcID}()$
	\State $[\M{Y}_{\lceil \log P \rceil}^{(p)},\M[\bar]{R}_{\lceil \log P \rceil}^{(p)}] = \text{Local-QR}(\M{A}^{(p)})$ \Comment{Leaf node QR}
  \If{$\lceil \log P \rceil \neq \lfloor \log P \rfloor$ } \Comment{Non-power-of-two case}
    \State $j = (p + 2^{\lfloor \log P \rfloor}) \mod 2^{\lceil \log P \rceil}$
    \If{$p \geq 2^{\lfloor \log P \rfloor}$} \Comment{Remainder processor}
      \State Send $\M[\bar]{R}_{\lceil \log P \rceil}^{(p)}$ to proc $j$
    \ElsIf{$p < P - 2^{\lfloor \log P \rfloor}$} \Comment{Partner of remainder processor}
      \State Receive $\M[\bar]{R}_{\lceil \log P \rceil}^{(j)}$ from proc $j$
      \State $[\M{Y}_{\star}^{(p)},\M[\bar]{R}_{\lceil \log P \rceil}^{(p)}] = \text{Local-QR}\left( \begin{bmatrix} \M[\bar]{R}_{\lceil \log P \rceil}^{(p)} \\ \M[\bar]{R}_{\lceil \log P \rceil}^{(j)} \end{bmatrix} \right)$
    \EndIf
  \EndIf
	\If{$p < 2^{\lfloor \log P\rfloor}$} \Comment{Butterfly tree on power-of-two procs}
    \For{$\ell= \lceil \log P \rceil -1$ down to $0$} 
	  	\State $j = 2^{\ell+1} \lfloor \frac{p}{2^{\ell+1}}\rfloor + \change{\left((p + 2^{\ell}) \mod 2^{\ell+1}\right)}$ \Comment{Determine partner}
	  	\State Send $\M[\bar]{R}_{\ell+1}^{(p)}$ to and receive $\M[\bar]{R}_{\ell+1}^{(j)}$ from proc $j$ \Comment{Communication}
	  	\If{$p <  j$}
	  		\State $[\M{Y}_{\ell}^{(p)},\M[\bar]{R}_{\ell}^{(p)}] = \text{Local-QR}\left( \begin{bmatrix} \M[\bar]{R}_{\ell+1}^{(p)} \\ \M[\bar]{R}_{\ell+1}^{(j)} \end{bmatrix} \right)$ \Comment{Tree node QR}
	  	\Else
	  		\State $[\M{Y}_{\ell}^{(p)},\M[\bar]{R}_{\ell}^{(p)}] = \text{Local-QR}\left( \begin{bmatrix} \M[\bar]{R}_{\ell+1}^{(j)} \\ \M[\bar]{R}_{\ell+1}^{(p)} \end{bmatrix} \right)$  \Comment{Partner tree node QR}
	  	\EndIf
	  \EndFor 
		\State $\M{R}=\bar{\M{R}}_0^{(p)}$
	\EndIf
	\If{$\lfloor \log P \rfloor \neq \lceil \log P \rceil$} \Comment{Non-power-of-two case}
		\State $j=(p + 2^{\lfloor \log P \rfloor}) \mod 2^{\lceil \log P \rceil}$ 
    \If{$p < P - 2^{\lfloor \log P \rfloor}$} \Comment{Partner of remainder proc}
			\State Send $\M{R}$ to proc $j$
    \ElsIf{$p \geq 2^{\lfloor \log P \rfloor}$}  \Comment{Remainder proc}
			\State Receive $\M{R}$ from proc $j$
		\EndIf
	\EndIf
\EndFunction
\end{algorithmic}
\end{algorithm}

\begin{algorithm}
\caption{Parallel Application of Implicit $Q$ from Butterfly TSQR}
\label{alg:BFappQapp}
\begin{algorithmic}[1]
\Require $\{\M{Y}_{\ell}^{(p)}\}$ represents \change{orthonormal} matrix $\M{Q}$ computed by \cref{alg:BFTSQRapp}
\Require $\M{C}$ is $b\times c$ and redundantly owned by all processors

\Ensure $\M{B}=\M{Q}\begin{bmatrix} \M{C} \\ \M{0} \end{bmatrix}$ is $m\times c$ and 1D-distributed so that proc $p$ owns row block $\M{B}^{(p)}$
\Function{$\M{B} = $ Par-TSQR-Apply-Q}{$\{\M{Y}_{\ell}^{(p)}\},\M{C}$}
\State $p = \textsc{MyProcID}()$
\If{$p < 2^{\lfloor \log P \rfloor}$} \Comment{Butterfly apply on power-of-two procs}
  \State $\bar{\M{B}}_0^{(p)} = \M{C}$
  \For{$\ell= 0$ to $ \lceil \log P \rceil -1$}
    \State $j = 2^{\ell+1} \lfloor \frac{p}{2^{\ell+1}}\rfloor + \change{\left((p + 2^{\ell}) \mod 2^{\ell+1}\right)}$ \Comment{Determine partner}
    \If{$p < j$}
      \State $\begin{bmatrix} \bar{\M{B}}_{\ell+1}^{(p)} \\ \M[\bar]{B}_{\ell+1}^{(j)} \end{bmatrix} = \textsc{Loc-Apply-Q}\left(\begin{bmatrix} \M{I}_b \\ \M{Y}_{\ell}^{(p)}\end{bmatrix}, \begin{bmatrix} \M[\bar]{B}_\ell^{(p)} \\ \M{0} \end{bmatrix} \right)$ \Comment{Tree node apply}
    \Else
      \State $\begin{bmatrix} \M[\bar]{B}_{\ell+1}^{(j)} \\ \M[\bar]{B}_{\ell+1}^{(p)} \end{bmatrix} = \textsc{Loc-Apply-Q}\left(\begin{bmatrix} \M{I}_b \\ \M{Y}_{\ell}^{(p)}\end{bmatrix}, \begin{bmatrix} \M[\bar]{B}_\ell^{(p)} \\ \M{0} \end{bmatrix} \right)$ \Comment{Partner apply}
    \EndIf
  \EndFor
  \EndIf
  \If{$\lfloor \log P \rfloor \neq \lceil \log P \rceil$} \Comment{Non-power-of-two case}
    \State $j=(p + 2^{\lfloor \log P \rfloor}) \mod 2^{\lceil \log P \rceil}$ 
    \If{$p < P - 2^{\lfloor \log P \rfloor}$} \Comment{Partner of remainder proc}
      \State $\begin{bmatrix} \bar{\M{B}}_{\lceil \log P \rceil}^{(p)} \\ \M[\bar]{B}_{\lceil \log P \rceil}^{(j)} \end{bmatrix} = \textsc{Loc-Apply-Q}\left(\begin{bmatrix} \M{I}_b \\ \M{Y}_{\star}^{(p)}\end{bmatrix}, \begin{bmatrix} \M[\bar]{B}_{\lceil \log P \rceil}^{(p)} \\ \M{0} \end{bmatrix} \right)$

      \State Send $\M[\bar]{B}_{\lceil \log P \rceil}^{(j)}$ to proc $j$
    \ElsIf{$p \geq 2^{\lfloor \log P \rfloor}$}  \Comment{Remainder proc}
      \State Receive $\M[\bar]{B}_{\lceil \log P \rceil}^{(p)}$ from proc $j$
    \EndIf
  \EndIf
  \State $\M{B}^{(p)} = \textsc{Loc-Apply-Q}\left(\M{Y}_{\lceil \log P \rceil}^{(p)}, \begin{bmatrix} \M[\bar]{B}_{\lceil \log P \rceil}^{(p)} \\ \M{0} \end{bmatrix} \right)$ \Comment{Leaf node apply}
\EndFunction
\end{algorithmic}
\end{algorithm}

\bibliography{main}
\bibliographystyle{siam}

\end{document}